# An oscillation-free fully partitioned scheme for the numerical modeling of cardiac active mechanics

F. Regazzoni[1] and A. Quarteroni[1,2]

[1]MOX - Dipartimento di Matematica, Politecnico di Milano,
P.zza Leonardo da Vinci 32, 20133 Milano, Italy
[2]Mathematics Institute, École Polytechnique Fédérale de Lausanne,
Av. Piccard, CH-1015 Lausanne, Switzerland (*Professor Emeritus*)

**Abstract**

In silico models of cardiac electromechanics couple together mathematical models describing different physics. One instance is represented by the model describing the generation of active force, coupled with the one of tissue mechanics. For the numerical solution of the coupled model, partitioned schemes, that foresee the sequential solution of the two subproblems, are often used. However, this approach may be unstable. For this reason, the coupled model is commonly solved as a unique system using Newton type algorithms, at the price, however, of high computational costs. In light of this motivation, in this paper we propose a new numerical scheme, that is numerically stable and accurate, yet within a fully partitioned (i.e. segregated) framework. Specifically, we introduce, with respect to standard segregated scheme, a numerically consistent stabilization term, capable of removing the nonphysical oscillations otherwise present in the numerical solution of the commonly used segregated scheme. Our new method is derived moving from a physics-based analysis on the microscale energetics of the force generation dynamics. By considering a model problem of active mechanics we prove that the proposed scheme is unconditionally absolutely stable (i.e. it is stable for any time step size), unlike the standard segregated scheme, and we also provide an interpretation of the scheme as a fractional step method. We show, by means of several numerical tests, that the proposed stabilization term successfully removes the nonphysical numerical oscillations characterizing the non stabilized segregated scheme solution. Our numerical tests are carried out for several force generation models available in the literature, namely the Niederer-Hunter-Smith model, the model by Land and coworkers, and the mean-field force generation model that we have recently proposed. Finally, we apply the proposed scheme in the context of a three-dimensional multiscale electromechanical simulation of the left ventricle.

# 1 Introduction

The heart can be regarded as a complex physiological system where different phenomena – involving different kind of physics: electrophysiology, biochemistry, mechanics, fluid-dynamics – interact at the organ, tissue and cellular scale [24, 25, 49]. One of the main challenges towards the construction of a multiphysics numerical model of the heart consists in coupling together the mathematical and numerical models describing the single physics involved in the cardiac function [6, 12, 32, 39, 47].

The currently available cardiac electromechanical models fall into two categories, namely loosely coupled and strongly coupled models [31, 37]. In loosely coupled models, only a one-way coupling between the different submodels is envisaged. In this case, the electrophysiological submodel is solved in a fixed domain. Then, the obtained calcium signal is used as an input for the submodel describing the generation of active force within cardiomyocytes. Finally, the displacement of the muscle is computed by using the resulting active force as an input. Despite being attractive – for their relatively small computational cost and since they allow to easily integrate stand-alone software tools already available for the different submodels – loosely coupled models neglect a complex system of feedback loops linking the different submodels in a nontrivial manner [4, 37]. As a matter of fact, the stretch of the cardiac tissue significantly affects the propagation of the electric potential through a series of phenomena collectively denoted as mechano-electrical feedback [29, 50]. Moreover, cardiac muscle cells are sensitive to mechanical stretch: the force generation phenomenon features several self-regulatory mechanisms associated with the feedback of the tissue strain (linked to the shortening of fibers) and of the strain rate (linked to the shortening velocity of fibers) [4, 17, 26]. In conclusion, while in some contexts loosely coupled models feature a sufficient level of accuracy, for many applications – that is when the feedback loops play a non negligible role – the submodels must be bidirectionally coupled, leading to the so-called strongly coupled models [6, 15, 39].

However, the construction of strongly coupled cardiac models clearly raises several challenges at the numerical level. The approaches to numerically couple different submodels can be classified into monolithic and segregated (or partitioned) schemes [31, 37]. Within a monolithic approach, one attempts to solve the entire coupled system at each time step simultaneously, typically by a Newton iterative scheme. Conversely, in segregated schemes, the different submodels are sequentially solved at each time step. In this manner, feedback information is only exchanged at the next time step. For this reason, monolithic approaches are more stable than the segregated ones [34, 37]. However, they feature a much higher computational cost, a more demanding implementation and do not allow to reuse pre-existing single-physics numerical codes in a straightforward manner. Alternatively, the submodels can be solved sequentially, but reiterating until a convergence criterion is satisfied before passing to the next time step [37]. This approach is however computationally demanding, as it may require several sub-iterations before reaching convergence. A further and significant advantage of segregated schemes is the possibility of using different spatial and temporal resolutions for the different submodels, in compliance with the characteristic spatial and temporal scales of each physics. This yields provide a great advantage from a computational standpoint with respect to monolithic schemes [11, 14, 39].

The main disadvantage of segregated schemes is that they may be unstable, especially when the role of the feedback – neglected at the numerical level during the solution of a single time step – becomes dominant. As a matter of fact, nonphysical oscillations typically show up when the force generation model and the mechanical model are treated sequentially (i.e. in a segregated manner) [31, 33, 34, 45, 56]. These instabilities are mainly linked to the feedback of the strain rate (i.e. the shortening velocity of the tissue) on the force generation model, as we will prove later. Due to the presence of these nonphysical oscillations, the force generation and the mechanical submodels are generally coupled in a monolithic way [29, 33, 34, 56]. As an alternative, in [31], the authors proposed to update only some of the terms of the force generation model after each Newton iteration of the mechanical model. This method, however, is not straightforward to implement (it is intrusive in the mechanics solver and cannot be applied to preexisting black-box solvers) and, most importantly, it is tailored on the model considered in [31] and cannot be easily generalized to different models, besides the ones written within the same formalism (the so-called *fading-memory model*).

In this paper, we mathematically analyze the source of these velocity-related instabilities, by considering a model with the minimal ingredients necessary to produce such numerical oscillations. Then, we propose a physically-inspired stabilized numerical scheme, still keeping a segregated approach between the force generation and the mechanical submodels.

## 1.1 Outline

This paper is structured as follows. First, in Sec. 2, we introduce the notation used in the paper. Then, in Sec. 3, we analyze the source of instabilities linked to the feedback from the tissue strain rate to the force generation model. Then, moving from energetic considerations, we derive the stabilized scheme that is proposed in this paper. In Sec. 4, we analyze the proposed scheme, by proving its numerical consistency and by providing a physical interpretation, an algebraic interpretation (thus proving its unconditional absolute stability) and, finally, an interpretation as a fractional step method. In Sec. 5 we show how the proposed scheme can be easily generalized to other force generation models than the one used for its derivation, and we test its effectiveness on three popular models available in the literature and on different test cases. These include a three-dimensional numerical simulation of the left ventricle. Moreover, we carry out a convergence analysis with respect to the time step size $\Delta t$. Finally, we draw our conclusions in Sec. 6.

## 2 Mathematical models and notation

Mathematical models of active force generation are typically written in the form of a system of ODEs. Specifically, denoting by $t$ the time variable, we consider force generation models written in the following form:

$$\begin{cases} \dot{\mathbf{r}}(t) = \mathbf{h}\left(\mathbf{r}(t), [\text{Ca}^{2+}]_i(t), \lambda(t), \dot{\lambda}(t)\right), & t \in (0, T], \\ T_a(t) = g\left(\mathbf{r}(t)\right), & t \in (0, T], \\ \mathbf{r}(0) = \mathbf{r}_0, \end{cases} \quad (1)$$

where $\mathbf{r}(t)$ is a vector collecting the state variables, which typically describe the state of the contractile proteins (several examples are provided in Appendices C, D and E). The inputs of the model are the intracellular calcium concentration $[\text{Ca}^{2+}]_i(t)$, the tissue strain in the fibers direction $\lambda(t)$ and its time derivative $\dot{\lambda}(t)$. The output of the model is $T_a(t)$, the active tension generated by the muscle tissue. Clearly, the functions $\mathbf{h}$ and $g$ depend on the model at hand. Where it is necessary for better clarity, we also use the notation $\frac{d}{dt}$ to denote time derivatives.

In the context of multiscale cardiac electromechanics, Eq. (1) is solved virtually in any point $\mathbf{x}$ of a computational domain $\Omega_0$, denoting the region occupied by the cardiac muscle tissue at rest. In practice, Eq. (1) is solved at each discretization node of the computational mesh or at each quadrature node. Then, Eq. (1) is coupled – on one side – with a model describing cardiac electrophysiology (e.g. the monodomain or bidomain equations [8]), which provides the input $[\text{Ca}^{2+}]_i(t)$, and – on the other side – with a model describing cardiac mechanics, written as

$$\begin{cases} \rho \dfrac{\partial^2 \mathbf{d}}{\partial t^2} - \nabla_0 \cdot (\mathbf{P}_{\text{pass}} + \mathbf{P}_{\text{act}}) = \mathbf{0}, & \text{in } \Omega_0 \times (0, T], \\ \textit{boundary conditions}, & \text{on } \partial\Omega_0 \times (0, T], \\ \mathbf{d} = \mathbf{d}_0, \quad \dfrac{\partial \mathbf{d}}{\partial t} = \mathbf{0}, & \text{in } \Omega_0 \times \{0\}, \end{cases} \qquad (2)$$

where $\Omega_0$ denotes the region occupied by the cardiac muscle tissue at rest, $\mathbf{d}\colon \Omega_0 \times [0, T] \to \mathbb{R}^3$ denotes its displacement and $\nabla_0 \cdot$ is the divergence operator in the reference domain $\Omega_0$. The Piola-Kirchhoff stress tensor is given by the sum of a passive term and an active term,

$$\mathbf{P}_{\text{pass}} = \frac{\partial \mathcal{W}}{\partial \mathbf{F}}, \qquad \mathbf{P}_{\text{act}} = T_\text{a} \frac{\mathbf{F}\mathbf{f}_0 \otimes \mathbf{f}_0}{|\mathbf{F}\mathbf{f}_0|}, \qquad (3)$$

where the strain energy density $\mathcal{W}\colon \text{Lin}^+ \to \mathbb{R}$ determines the passive constitutive behavior of the tissue, $\mathbf{f}_0$ denotes the direction of cardiac fibers and where the magnitude of active tension $T_\text{a}$ is provided by the force generation model (1). In some cases, with a quasistatic approximation, the inertia term is neglected (i.e. $\rho = 0$). In other cases, viscous terms are introduced either in the constitutive law (viscoelastic models) or in the boundary conditions [39, 58].

The solution of model (2) allows to compute the strain in the fibers direction, defined as

$$\lambda(t) = \lambda(\mathbf{d}(t)) = \sqrt{\mathcal{I}_{4,f}(\mathbf{F}(t))} - 1, \qquad (4)$$

where $\mathcal{I}_{4,f}(\mathbf{F}) = \mathbf{F}\mathbf{f}_0 \cdot \mathbf{F}\mathbf{f}_0$. We notice that we have $\lambda(t) = 0$ when the tissue is at rest, $\lambda(t) > 0$ when it is stretched, $\lambda(t) < 0$ when it is compressed. This provides a feedback to the force generation model (1), whose dynamics depends on both $\lambda(t)$ and $\dot{\lambda}(t)$.

When we consider time discretization, we denote by $\Delta t$ the time step length and we define a discrete collection of times as $t_k = k\,\Delta t$, for $k = 0, 1, 2, \ldots$. Then, we denote with a superscript $(k)$ the approximation of a given variable at time $t_k$. For instance, we have $\mathbf{r}^{(k)} \simeq \mathbf{r}(t_k)$ and $T_\text{a}^{(k)} \simeq T_\text{a}(t_k)$.

Within the segregated approach, at each time step $k$ we first update the solution of the force generation model, by considering, for instance, the following implicit Euler scheme:

$$\begin{cases} \dfrac{\mathbf{r}^{(k+1)} - \mathbf{r}^{(k)}}{\Delta t} = \mathbf{h}\left(\mathbf{r}^{(k+1)}, [\text{Ca}^{2+}]_\text{i}^{(k+1)}, \lambda^{(*)}, \dot{\lambda}^{(*)}\right), & k \geq 0, \\ T_\text{a}^{(k+1)} = g\left(\mathbf{r}^{(k+1)}\right), & k \geq 0, \\ \mathbf{r}^{(0)} = \mathbf{r}_0 \end{cases} \qquad (5)$$

where we employ, as an approximation of $\lambda(t_k)$ and $\dot{\lambda}(t_k)$, a first-order extrapolation from the previous time steps (higher-order extrapolations can be considered as well), by setting

$$\lambda^{(*)} = \lambda^{(k)}, \quad \dot{\lambda}^{(*)} = \frac{\lambda^{(k)} - \lambda^{(k-1)}}{\Delta t},$$

where we denote $\lambda^{(k)} = \sqrt{\mathcal{I}_{4,f}(\mathbf{F}^{(k)})} - 1$. Then, we employ the updated value of the active tension to update the displacement variable, by means – for instance – of the following BDF (backward differentiation formula) implicit scheme:

$$\begin{cases} \rho \dfrac{\mathbf{d}^{(k+1)} - 2\mathbf{d}^{(k)} + \mathbf{d}^{(k-1)}}{\Delta t^2} - \nabla_0 \cdot \left( \mathbf{P}_{\text{pass}}(\mathbf{F}^{(k+1)}) + \mathbf{P}_{\text{act}}^{(k+1)} \right) = \mathbf{0}, \qquad k \geq 0, \\ \mathbf{d}^{(0)} = \mathbf{d}^{(-1)} = \mathbf{d}_0, \end{cases} \quad (6)$$

where

$$\mathbf{P}_{\text{act}}^{(k+1)} = T_{\text{a}}^{(k+1)} \frac{\mathbf{F}^{(k+1)} \mathbf{f}_0 \otimes \mathbf{f}_0}{|\mathbf{F}^{(k+1)} \mathbf{f}_0|}, \quad (7)$$

and where space discretization is addressed with a suitable method, such as the Finite Element Method (FEM), also accounting for the boundary conditions.

Conversely, within the monolithic approach, the mechanical problem (6) is solved simultaneously with the force generation model, by setting in (5)

$$\lambda^{(*)} = \lambda^{(k+1)}, \quad \dot{\lambda}^{(*)} = \frac{\lambda^{(k+1)} - \lambda^{(k)}}{\Delta t}.$$

## 3 Analysis of the source of instabilities

With the aim of gaining physical insight into the source of velocity-dependent instabilities linked to the feedback of mechanical deformation on the force generation dynamics, we look for a minimal model that contains all the necessary ingredients to generate such instabilities.

### 3.1 The Huxley 1957 model (H57)

When a muscle fiber contracts at a constant length (i.e. isometrically), it generates a larger force than when its length is decreasing. In other terms, the generated force ($T_{\text{a}}$) is a decreasing function of the shortening velocity ($-\dot{\lambda}$) [18, 49]. The mechanisms underlying this phenomenon have been firstly revealed and explained in the celebrated Huxley model [23] (denoted henceforth as H57 model).

#### 3.1.1 Derivaton of the H57 model

Huxley considered a population of myosin heads and actin binding sites, large enough to assume that the probability of finding an actin-myosin pair at distance $x$ is constant in an interval sufficiently close to $x = 0$. More precisely, we denote by $\rho_{\text{AM}}$ the linear density of actin-myosin pairs in each pair of interacting thin and thick filaments. This variable that can be computed as $\rho_{\text{AM}} = L_{\text{so}}/(D_M D_A)$, where $L_{\text{so}}$ is the single-overlap length (i.e. the length of the region of overlap between the thin and the thick filaments), while $D_M$ (respectively, $D_A$) represents the distance between two consecutive myosin heads (respectively, actin binding sites).

Each actin-myosin pair has two possible states: when detached they do not interact; when attached, they form a crossbridge (XB) that behaves as a linear spring, generating a force equal to $k_{\text{XB}} x$ (by convention, $x$ is positive when attachment leads to a positive force). The attachment and detachment rates are strain-dependent and are assigned by the functions $f(x)$ and $g(x)$, respectively. Clearly, $f(x)$ is large for $x$ close to 0 and it vanishes far from $x = 0$, while $g(x)$ has an opposite behavior.

We denote by $n(x,t) \in [0,1]$ the probability that an actin-myosin pair with distance $x$ is attached. Hence, the total number – in half sarcomere – of attached XBs with elongation between $a$ and $b$ is $\int_a^b \rho_{\text{AM}} n(x,t)\,dx$. The probability density $n(x,t)$ is convected by the mutual sliding velocity between the thin and the thick filament, given by $v_{\text{hs}}(t) = -\frac{dSL(t)/2}{dt}$, where $SL(t)$ denotes the current sarcomere length. In conclusion, it can be proved (see e.g. [45]) that the evolution of $n(x,t)$ satisfies the following PDE (H57 model):

$$\frac{\partial n(x,t)}{\partial t} - v_{\text{hs}}(t)\frac{\partial n(x,t)}{\partial x} = (1 - n(x,t))f(x) - n(x,t)g(x), \qquad x \in \mathbb{R},\, t \geq 0, \quad (8)$$

endowed with a suitable initial condition. Denoting by $\sigma_{\text{hf}}$ the area density of pairs of interacting thin filaments and thick filaments, the tissue-level active tension $T_a(t)$ can be obtained as the product of $\sigma_{\text{hf}}$ times the force generated by each pair of filaments. Therefore, we have

$$T_a(t) = \sigma_{\text{hf}}\,\rho_{\text{AM}}\,k_{\text{XB}} \int_{-\infty}^{+\infty} x\, n(x,t)\,dx. \quad (9)$$

### 3.1.2 Distribution-moments equations

Under suitable hypotheses on the transition rate functions $f(x)$ and $g(x)$, the PDE (8) can be reduced to an ODE [5, 7, 57]. In fact, let us introduce the dimensionless distribution-moments of the density function $n(x,t)$ (for $p = 0, 1, \ldots$):

$$\mu_p(t) := \int_{-\infty}^{+\infty} \left(\frac{x}{SL_0/2}\right)^p n(x,t)\frac{dx}{D_M}, \quad (10)$$

and of the transition rate function $f(x)$:

$$\mu_p^f := \int_{-\infty}^{+\infty} \left(\frac{x}{SL_0/2}\right)^p f(x) \frac{dx}{D_M}.$$

Under the (physically motivated [7]) hypothesis that $f(x) + g(x) = r$ (a constant) for any $x$, the PDE (8) reduces to the following pair of ODEs:

$$\begin{cases} \dot{\mu}_0(t) = \mu_0^f - r\,\mu_0(t), & t \geq 0, \\ \dot{\mu}_1(t) = \mu_1^f - r\,\mu_1(t) + \dot{\lambda}(t)\,\mu_0(t), & t \geq 0, \\ \mu_0(0) = \mu_{0,0}, \quad \mu_1(0) = \mu_{1,0}, \end{cases} \quad (11)$$

where we have used the fact that $SL = SL_0(1 + \lambda)$, being $\lambda$ the strain in the fibers direction and $SL_0$ the slack sarcomere length. The constant $\mu_{0,0}$ and $\mu_{1,0}$ denote the initial values of the two variables. For physical meaningfulness, we have $\mu_{0,0} \in [0, \mu_0^f/r]$, which clearly entails $\mu_0(t) \in [0, \mu_0^f/r]$ for any $t \geq 0$. From Eq. (9) it follows that

$$T_a(t) = a_{\text{XB}}\,\mu_1(t), \quad (12)$$

where $a_{\text{XB}} = \sigma_{\text{hf}}\,\rho_{\text{AM}}\,D_M\,k_{\text{XB}}\,\frac{SL_0}{2}$ represents the XB stiffness upscaled at the tissue level. The zero-order moment $\mu_0(t)$ can be interpreted as the fraction of actin binding sites that are involved in a XB, while $\mu_1(t)/\mu_0(t)$ represents the mean elongation of attached XBs. This allows a better physical understanding of Eq. (12).

The distribution-moments model (11) can be considered as a minimal model for the feedback of shortening velocity on force generation. In fact, let us suppose to start from the equilibrium configuration $\mu_0(0) = \mu_0^f/r$, $\mu_1(0) = \mu_1^f/r$ (isometric conditions) and let us shorten the tissue with a constant velocity $\dot{\lambda} < 0$. This clearly decreases the elongation of attached XB, thus lowering the first-order moment $\mu_1$. In fact, the steady-state solution is given by $\mu_1 = \mu_1^f/r + \mu_0^f/r^2 \dot{\lambda} < \mu_1(0)$. In conclusion, Eq. (12) entails that the steady-state active tension decreases with the shortening velocity of the muscle fibers, coherently with what observed by A. V. Hill [18].

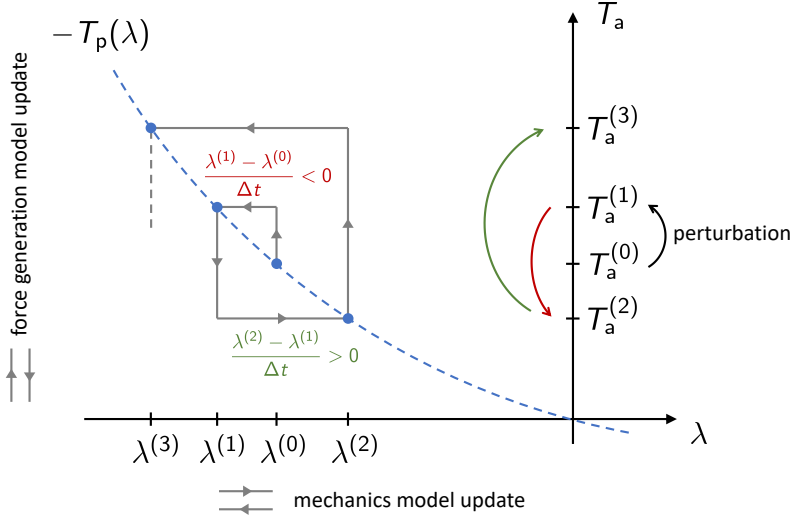

Figure 1: Visual representation of the feedback loop that generates nonphysical oscillations when a force generation model and a passive mechanics model are coupled in a segregated manner. Denoting by $T_{\mathrm{p}}(\lambda)$ the passive tension, the equilibrium configuration is obtained for $T_{\mathrm{a}}(t) = -T_{\mathrm{p}}(\lambda(t))$ (the latter curve is represented by the blue dashed line). Within a the segregated approach, at each time step one first updates the force generation model (by moving along the vertical axis) and then updates the mechanical model (by moving along the horizontal axis).

Moreover, as we will numerically show in Sec. 4.3, the distribution-moments model (11) produces nonphysical numerical oscillations when coupled in a segregated manner even with a simple 0D model of tissue mechanics. The mechanism underlying these velocity-related instabilities is rooted in the feedback loop schematically illustrated in Fig. 1. Suppose that at the iteration $k = 0$ the tissue is at rest, with strain $\lambda^{(0)}$ and active tension $T_{\mathrm{a}}^{(0)}$. Suppose then that a small perturbation (due e.g. to the excitation of the tissue) makes the active tension increase ($T_{\mathrm{a}}^{(1)} > T_{\mathrm{a}}^{(0)}$) and, as a consequence, the fibers shorten ($\lambda^{(1)} < \lambda^{(0)}$). In the next time step, the negative strain rate ($\frac{\lambda^{(1)} - \lambda^{(0)}}{\Delta t} < 0$) causes a drop of active force ($T_{\mathrm{a}}^{(2)} < T_{\mathrm{a}}^{(1)}$) and, consequently, of the shortening of the fiber ($\lambda^{(2)} > \lambda^{(1)}$). In the following time step, because of the positive strain rate, the active force raises again ($T_{\mathrm{a}}^{(3)} > T_{\mathrm{a}}^{(2)}$) and the cycle is repeated.

### 3.2 Energetic analysis of the segregated scheme

In this section, we examine the source of the above mentioned velocity-related instabilities for the model (11), based on energetic considerations.

### 3.2.1 Deriving the active Piola-Kirchhoff stress tensor from microscale energetics

Going back to the microscopic derivation of the model (11), we recall that each attached XB generates a force equal to $k_{\text{XB}} x$ (see Sec. 3.1). Hence, we can associate to each attached XB the following elastic energy:

$$\mathcal{W}_{\text{XB}}(x) = \frac{1}{2} k_{\text{XB}} x^2.$$

Therefore, the total elastic energy associated with a pair of interacting thin and thick filaments is given by

$$\mathcal{W}_{\text{hf}}(t) = \frac{1}{2} \rho_{\text{AM}}\, k_{\text{XB}} \int_{-\infty}^{+\infty} x^2\, n(x,t)\, dx.$$

Due to the additive nature of energy, the total tissue-level energy density of attached XBs (i.e. energy per unit volume) is given by the product of $\mathcal{W}_{\text{hf}}(t)$ times the volume density of interacting filaments, given by $\rho_{\text{hf}} = \sigma_{\text{hf}}/(SL_0/2)$:

$$\begin{aligned}
\mathcal{W}_{\text{act}}(t) &= \frac{1}{SL_0} \sigma_{\text{hf}}\, \rho_{\text{AM}}\, k_{\text{XB}} \int_{-\infty}^{+\infty} x^2\, n(x,t)\, dx \\
&= \frac{1}{4} \sigma_{\text{hf}}\, \rho_{\text{AM}}\, k_{\text{XB}}\, SL_0\, D_M \int_{-\infty}^{+\infty} \left(\frac{x}{SL_0/2}\right)^2 n(x,t)\, \frac{dx}{D_M} \\
&= \frac{1}{2} a_{\text{XB}}\, \mu_2(t),
\end{aligned}$$

where we have used the definition of Eq. (10).

We recall from Eq. (3) that the Piola-Kirchhoff stress tensor of a deformed elastic body is obtained by differentiating the elastic energy with respect to the strain tensor $\mathbf{F}$. Similarly, since we model attached XBs as elastic elements, we can compute the active stress generated by them by differentiating $\mathcal{W}_{\text{act}}(t)$ with respect to $\mathbf{F}$. However, for this operation to be meaningful, it must be carried out in a strictly Lagrangian formalism (i.e. referred to the reference configuration). Nonetheless, the variable $x$ (used to define the moments $\mu_p$) is rather an Eulerian variable (i.e. it is referred to the current configuration). Hence, we consider its pullback in the reference configuration:

$$\hat{x} = x - \frac{SL - SL_0}{2} = x - \frac{SL_0}{2} \lambda,$$

and we define the probability density function in Lagrangian coordinates as

$$\hat{n}(\hat{x}, t) = n\left(\hat{x} + \frac{SL_0}{2} \lambda, t\right).$$

Moreover, we define the moments of the displacement in Lagrangian coordinates, for $p \in \mathbb{N}$ as

$$\begin{aligned}
\hat{\mu}_p(t) &= \int_{-\infty}^{+\infty} \left(\frac{\hat{x}}{SL_0/2}\right)^p \hat{n}(\hat{x}, t)\, \frac{d\hat{x}}{D_M} \\
&= \int_{-\infty}^{+\infty} \left(\frac{x}{SL_0/2} - \lambda\right)^p n(x,t)\, \frac{dx}{D_M} \\
&= \sum_{k=0}^{p} \binom{p}{k} (-\lambda)^{p-k} \mu_k(t) \\
&= \mu_p(t) - p\,\lambda\,\mu_{p-1}(t) + \cdots + p\,(-\lambda)^{p-1}\,\mu_1(t) + (-\lambda)^p\,\mu_0(t).
\end{aligned}$$

It follows that the moments $\hat{\mu}_p$ and $\mu_p$ are linked by the following relationships:

$$\mu_0 = \hat{\mu}_0, \qquad \mu_1 = \lambda\hat{\mu}_0 + \hat{\mu}_1, \qquad \mu_2 = \lambda^2\hat{\mu}_0 + 2\lambda\hat{\mu}_1 + \hat{\mu}_2. \tag{13}$$

This entails

$$\mathcal{W}_{\text{act}} = \frac{1}{2} a_{\text{XB}} \left( \hat{\mu}_0 \lambda^2 + 2\hat{\mu}_1 \lambda + \hat{\mu}_2 \right). \tag{14}$$

Equation (14) provides an expression for the energy associated with attached XBs in a fully Lagrangian frame of reference, in which the unique term depending on $\mathbf{F}$ is $\lambda = \lambda(\mathbf{F})$. Therefore, the active part of the Piola-Kirchhoff stress tensor can be obtained by the chain rule as:

$$\mathbf{P}_{\text{act}} = \frac{\partial \mathcal{W}_{\text{act}}}{\partial \mathbf{F}} = \frac{\partial \mathcal{W}_{\text{act}}}{\partial \lambda} \frac{\partial \lambda}{\partial \mathbf{F}},$$

where:

$$\frac{\partial \mathcal{W}_{\text{act}}}{\partial \lambda} = a_{\text{XB}} \left( \hat{\mu}_0 \lambda + \hat{\mu}_1 \right), \qquad \frac{\partial \lambda}{\partial \mathbf{F}} = \frac{\mathbf{F}\mathbf{f}_0 \otimes \mathbf{f}_0}{1+\lambda}.$$

Finally, we obtain:

$$\mathbf{P}_{\text{act}} = a_{\text{XB}} \left( \hat{\mu}_0 \lambda + \hat{\mu}_1 \right) \frac{\mathbf{F}\mathbf{f}_0 \otimes \mathbf{f}_0}{|\mathbf{F}\mathbf{f}_0|}. \tag{15}$$

By comparing Eq. (15) with Eq. (3), we obtain $T_{\text{a}} = a_{\text{XB}}(\hat{\mu}_0 \lambda + \hat{\mu}_1)$ and thus (by Eq. (13)) $T_{\text{a}} = a_{\text{XB}}\mu_1$, coherently with Eq. (12). We have thus obtained an equivalent derivation of Eq. (12), uniquely based on energetic considerations.

In conclusion, we have two equivalent formulations for the active part of the Piola-Kirchhoff stress tensor:

(F1) $\mathbf{P}_{\text{act}} = a_{\text{XB}} \, \mu_1 \, \dfrac{\mathbf{F}\mathbf{f}_0 \otimes \mathbf{f}_0}{|\mathbf{F}\mathbf{f}_0|}$;

(F2) $\mathbf{P}_{\text{act}} = a_{\text{XB}} \left( \hat{\mu}_1 + \hat{\mu}_0 \lambda \right) \dfrac{\mathbf{F}\mathbf{f}_0 \otimes \mathbf{f}_0}{|\mathbf{F}\mathbf{f}_0|}$.

The difference between the two formulations is in the coordinate system (Lagrangian vs Eulerian, see e.g. [1]), used to describe the microscopic elongation of the myosin arms. Indeed, (F1) refers to an hybrid Lagrangian-Eulerian formalism: while the macroscale strain $\mathbf{F}$ is written in Lagrangian coordinates, the variable $\mu_1$ is defined as the first-order distribution-moment of the microscopic Eulerian coordinate $x$. On the other hand, (F2) is a fully Lagrangian formalism (i.e. both the macroscopic strain $\mathbf{F}$ and the microscopic strain variable $\hat{x}$ are referred to the reference configuration).

### 3.2.2 A fully Lagrangian active mechanics scheme

In the previous section we have shown that the active Piola-Kirchhoff stress tensor of (F2) can be interpreted as the differential, with respect to the strain tensor, of the elastic energy associated with attached XBs. Formulation (F2) is clearly equivalent to formulation (F1), which is typically used in the literature. In the latter, however, the dependence of the active tension on the strain $\lambda$ is somehow hidden.

From Eq. (13), it follows that the Lagrangian moments can be derived from the eulerian ones as

$$\hat{\mu}_0(\lambda) = \mu_0, \quad \hat{\mu}_1(\lambda) = \mu_1 - \lambda\mu_0. \tag{16}$$

Within a segregated scheme, the values of $\mu_0^{(k+1)}$ and $\mu_1^{(k+1)}$ are obtained by employing the value of $\lambda$ at the iteration $k$. Hence, based on (16), we define the corresponding Lagrangian variables as

$$\hat{\mu}_0^{(k+1)} = \hat{\mu}_0(\lambda^{(k)}) = \mu_0^{(k+1)}, \quad \hat{\mu}_1^{(k+1)} = \hat{\mu}_1(\lambda^{(k)}) = \mu_1^{(k+1)} - \lambda^{(k)}\mu_0^{(k+1)}. \tag{17}$$

It follows that the standard segregated scheme of Eq. (7), in a fully Lagrangian frame of reference, reads as follows:

$$\mathbf{P}_{\text{act}}^{(k+1)} = a_{\text{XB}} \left(\hat{\mu}_1^{(k+1)} + \hat{\mu}_0^{(k+1)}\lambda^{(k)}\right) \frac{\mathbf{F}^{(k+1)}\mathbf{f}_0 \otimes \mathbf{f}_0}{|\mathbf{F}^{(k+1)}\mathbf{f}_0|}. \tag{18}$$

Hence, despite Eq. (7) is formally written as a fully-implicit scheme (that is typically unconditionally stable for problems of this kind [38]), it actually hides an implicit-explicit scheme, where the strain $\lambda$ is treated explicitly (in fact we have $\lambda^{(k)}$ rather than $\lambda^{(k+1)}$). Therefore, we consider the following discrete-in-time Piola-Kirchhoff stress tensor written in a fully Lagrangian formalism (i.e. based on formulation (F2)):

$$\begin{aligned}\mathbf{P}_{\text{act}}^{(k+1)} &= a_{\text{XB}} \left(\hat{\mu}_1^{(k+1)} + \hat{\mu}_0^{(k+1)}\lambda^{(k+1)}\right) \frac{\mathbf{F}^{(k+1)}\mathbf{f}_0 \otimes \mathbf{f}_0}{|\mathbf{F}^{(k+1)}\mathbf{f}_0|} \\ &= a_{\text{XB}} \left[\mu_1^{(k+1)} + \mu_0^{(k+1)}\left(\lambda^{(k+1)} - \lambda^{(k)}\right)\right] \frac{\mathbf{F}^{(k+1)}\mathbf{f}_0 \otimes \mathbf{f}_0}{|\mathbf{F}^{(k+1)}\mathbf{f}_0|},\end{aligned}$$

that is

$$\mathbf{P}_{\text{act}}^{(k+1)} = \left[T_{\text{a}}^{(k+1)} + K_{\text{a}}^{(k+1)}\left(\lambda^{(k+1)} - \lambda^{(k)}\right)\right] \frac{\mathbf{F}^{(k+1)}\mathbf{f}_0 \otimes \mathbf{f}_0}{|\mathbf{F}^{(k+1)}\mathbf{f}_0|}, \tag{19}$$

where $K_{\text{a}}^{(k+1)} = a_{\text{XB}}\mu_0^{(k+1)}$ is the *active stiffness* of the tissue. Indeed, since $\mu_0(t)$ represents the fraction of actin binding sites involved in a XB, the term $K_{\text{a}}(t) = a_{\text{XB}}\mu_0(t)$ represents the total stiffness (at the tissue level) of attached XBs.

In what follows, we will refer to the numerical scheme consisting of (5), (6) and (19) (in substitution of (7)) as the *stabilized-segregated scheme*, for reasons that will be clear later.

## 4 Analysis of the proposed stabilized scheme

In this section we analyze the stabilized-segregated scheme.

### 4.1 Numerical consistency of the stabilized-segregated scheme

The stabilized-segregated scheme is *numerically consistent* (in the sense of [36]) with the problem (1)–(2). Indeed, by setting the discretized variables equal to the exact solution (i.e. $\mathbf{d}_h^{(k)} = \mathbf{d}(t_k)$, $T_{\text{a}}^{(k)} = T_{\text{a}}(t_k)$ and $K_{\text{a}}^{(k)} = K_{\text{a}}(t_k)$) and by letting $\Delta t \to 0$, we get

$$\left[T_{\text{a}}^{(k)} + K_{\text{a}}^{(k)}\left(\lambda^{(k+1)} - \lambda^{(k)}\right)\right] \xrightarrow{\Delta t \to 0} T_{\text{a}}(t_k).$$

As a matter of fact, we have $\left(\lambda^{(k+1)} - \lambda^{(k)}\right) = \left(\mathcal{I}_{4,f}(\mathbf{F}^{(k+1)}) - \mathcal{I}_{4,f}(\mathbf{F}^{(k)})\right) = \mathcal{O}(\Delta t)$. The newly introduced term can thus be interpreted as a consistent stabilization term (of first order with respect to $\Delta t$).

## 4.2 Physical interpretation of the stabilized-segregated scheme

The stabilized-segregated scheme that we propose consists in replacing, in the momentum equation, the term $T_a^{(k+1)}$ by the term $[T_a^{(k+1)} + K_a^{(k+1)}(\lambda^{(k+1)} - \lambda^{(k)})]$. This means that in the stabilized-segregated scheme the active tension is not seen – at the numerical level – as a constant force. Rather, it is regarded as an elastic force, whose value depends on the tissue strain in the fibers direction. This is more coherent with the microscopical basis of the force generation model, in which each attached XB behaves as a linear spring (see Sec. 3.1.1).

We remark that formulation (19) does not introduce a new model with respect to (7). As a matter of fact, the difference is only at the numerical level (see Sec. 4.1), that is we are using a different, but still consistent, numerical scheme. Nonetheless, the stabilized-segregated scheme better reflects the physics underlying the model, thus featuring better numerical stability properties, as we prove later.

## 4.3 Algebraic interpretation of the stabilized-segregated scheme

In order to provide an algebraic interpretation of the newly introduces term, we couple the minimal model of force generation of Eq. (11) with a minimal model of tissue mechanics, in replacement of Eq. (2). Specifically, we consider the following zero-dimensional model for the tissue strain $\lambda$:

$$M\ddot{\lambda}(t) + \sigma\dot{\lambda}(t) + K_p\lambda(t) + T_a(t) = p(t), \quad t \geq 0, \tag{20}$$

with suitable initial conditions on $\lambda(0)$ and $\dot{\lambda}(0)$ and where $M$ represents a normalized mass, $\sigma$ a normalized viscous modulus, $K_p$ the passive stiffness of the tissue and $p(t)$ an externally applied load. In what follows, we will consider also the quasistatic approximation of Eq. (20) (i.e. by neglecting inertia and viscous damping):

$$K_p\lambda(t) + T_a(t) = p(t), \quad t \geq 0. \tag{21}$$

In conclusion, the minimal model of active mechanics that we will consider is

$$\begin{cases} \dot{\mu}_0(t) = \mu_0^f - r\,\mu_0(t) & t \geq 0, \\ \dot{\mu}_1(t) = \mu_1^f - r\,\mu_1(t) + \dot{\lambda}(t)\,\mu_0(t) & t \geq 0, \\ M\ddot{\lambda}(t) + \sigma\dot{\lambda}(t) + K_p\lambda(t) = p(t) - a_{\text{XB}}\,\mu_1(t) & t \geq 0, \\ \mu_0(0) = \mu_{0,0}, \quad \mu_1(0) = \mu_{1,0}, \quad \lambda(0) = \lambda_0, \quad \dot{\lambda}(0) = 0, \end{cases} \tag{22}$$

where we set, in the quasistatic case, $M = \sigma = 0$.

### 4.3.1 Quasistatic elasticity case

Let us first consider the quasistatic case (i.e. by setting $M = \sigma = 0$). We consider and compare the following three strategies for the time-discretization of Eq. (22) and for the coupling between the two submodels.

**Monolithic scheme.** Within the monolithic strategy, we simultaneously solve for the activation variables ($\mu_0$ and $\mu_1$) and for $\lambda$:

$$\begin{cases} \dfrac{\mu_0^{(k+1)} - \mu_0^{(k)}}{\Delta t} = \mu_0^f - r\,\mu_0^{(k+1)} & k = 0, 1, \ldots \\ \dfrac{\mu_1^{(k+1)} - \mu_1^{(k)}}{\Delta t} = \mu_1^f - r\,\mu_1^{(k+1)} + \dfrac{\lambda^{(k+1)} - \lambda^{(k)}}{\Delta t}\mu_0^{(k+1)} & k = 0, 1, \ldots \\ K_p \lambda^{(k+1)} = p(t_{k+1}) - a_{\mathrm{XB}}\,\mu_1^{(k+1)} & k = 0, 1, \ldots \\ \mu_0^{(0)} = \mu_{0,0}, \quad \mu_1^{(0)} = \mu_{1,0}, \quad \lambda^{(0)} = \lambda_0. \end{cases} \quad (23)$$

**Segregated scheme.** In the segregated scheme, we initialize the variables as $\mu_0^{(0)} = \mu_{0,0}$, $\mu_1^{(0)} = \mu_{1,0}$ and $\lambda^{(-1)} = \lambda^{(0)} = \lambda_0$. At each time step, we first approximate the solution of the force generation model by employing the values of $\lambda$ obtained at the previous time step:

$$\begin{cases} \dfrac{\mu_0^{(k+1)} - \mu_0^{(k)}}{\Delta t} = \mu_0^f - r\,\mu_0^{(k+1)} & k = 0, 1, \ldots \\ \dfrac{\mu_1^{(k+1)} - \mu_1^{(k)}}{\Delta t} = \mu_1^f - r\,\mu_1^{(k+1)} + \dfrac{\lambda^{(k)} - \lambda^{(k-1)}}{\Delta t}\mu_0^{(k+1)} & k = 0, 1, \ldots \end{cases} \quad (24)$$

Then, we update the value of $\lambda$ by approximating the mechanics model:

$$K_p \lambda^{(k+1)} = p(t_{k+1}) - a_{\mathrm{XB}}\,\mu_1^{(k+1)} \quad k = 0, 1, \ldots \quad (25)$$

**Stabilized-segregated scheme.** The initialization and the first step of the proposed stabilized scheme coincides with that of the segregated scheme (see Eq. (24)). However, the second step is modified as follows:

$$K_p \lambda^{(k+1)} = p(t_{k+1}) - a_{\mathrm{XB}} \left[ \mu_1^{(k+1)} + \mu_0^{(k+1)} \left( \lambda^{(k+1)} - \lambda^{(k)} \right) \right] \quad k = 0, 1, \ldots \quad (26)$$

Let us analyze the temporal stability of the three schemes. We notice that each of these schemes can be written in the following form:

$$\begin{cases} \boldsymbol{\psi}^{(k)} = \Phi(\boldsymbol{\psi}^{(k-1)}, t_k, \Delta t) & k = 1, 2, \ldots \\ \boldsymbol{\psi}^{(0)} = \boldsymbol{\psi}_0, \end{cases} \quad (27)$$

where we have defined the state vector $\boldsymbol{\psi}^{(k)} = (\mu_0^{(k)}, \mu_1^{(k)}, \lambda^{(k)}, \lambda^{(k-1)})^T$, with initial value $\boldsymbol{\psi}_0 = (\mu_{0,0}, \mu_{1,0}, \lambda_0, \lambda_0)^T$, and where $\Phi \colon \mathbb{R}^n \times \mathbb{R}^+ \times \mathbb{R}^+ \to \mathbb{R}^n$ denotes an iteration function across the time steps. It should not be surprising that the state $\boldsymbol{\psi}^{(k)}$ contains two consecutive values of $\lambda^{(k)}$, as the dynamics of $\mu_1$ depends on the time derivative of $\lambda$, approximated by an incremental quotient. In order to study the stability of the schemes written in the form (27), we introduce the perturbed problem

$$\begin{cases} \widetilde{\boldsymbol{\psi}}^{(k)} = \Phi(\widetilde{\boldsymbol{\psi}}^{(k-1)}, t_k, \Delta t) + \Delta t\,\boldsymbol{\eta}^{(k)} & k = 1, 2, \ldots \\ \widetilde{\boldsymbol{\psi}}^{(0)} = \boldsymbol{\psi}_0 + \boldsymbol{\eta}^{(0)}, \end{cases} \quad (28)$$

where $\boldsymbol{\eta}^{(k)}$ denotes a suitable perturbation. We recall the following definitions [36]:

**Definition 1** (zero-stability). *Let us consider a finite time $T > 0$, and let us consider a uniform subdivision of the time interval $(0, T)$ into $N$ subintervals, i.e. $\Delta t = T/N$. The numerical scheme (27) is zero-stable if*

$$\exists \, \Delta t_0 > 0, \, \exists \, C > 0 : \quad \sup_{k=0,\ldots,N} |\widetilde{\boldsymbol{\psi}}^{(k)} - \boldsymbol{\psi}^{(k)}| \leq C \sup_{k=0,\ldots,N} |\boldsymbol{\eta}^{(k)}| \quad \forall \, \Delta t \in (0, \Delta t_0],$$

*where $\boldsymbol{\psi}^{(k)}$ and $\widetilde{\boldsymbol{\psi}}^{(k)}$ are solutions to problems (27) and (28), respectively.*

**Definition 2** (absolute stability). *Let us consider a given time step size $\Delta t > 0$, and let us consider the solution of the numerical problem for $t_k \to +\infty$. The numerical scheme (27) is absolutely stable in correspondence to $\Delta t$ if*

$$\exists \, C > 0 : \quad \lim_{k \to +\infty} |\widetilde{\boldsymbol{\psi}}^{(k)} - \boldsymbol{\psi}^{(k)}| \leq C \sup_{k=0,\ldots} |\boldsymbol{\eta}^{(k)}|,$$

*where $\boldsymbol{\psi}^{(k)}$ and $\widetilde{\boldsymbol{\psi}}^{(k)}$ are solutions to problems (27) and (28), respectively. Moreover, we say that the method is unconditionally absolutely stable if it is absolutely stable for any $\Delta t > 0$.*

On the one hand, the zero-stability deals with the behavior of the numerical scheme on a given temporal interval for $\Delta t \to 0$ (whence the name of *zero*-stability). This request arises from the need of keeping under control the unavoidable round-off errors due to the finite arithmetic of computers (represented in this context by the perturbations $\boldsymbol{\eta}^{(k)}$). On the other hand, the absolute stability deals with the behavior of the method for a given $\Delta t$ in the limit $t_k \to +\infty$. This property guarantees that the effect of perturbations is kept under control and does not lead to a blow-up of the solution. For a dissipative (Lyapunov asymptotically stable) dynamical system, absolute stability guarantees that the numerical solution will tend to zero as $t_k \to +\infty$ if the only perturbation concerns the initial datum [36].

The following result (whose proof is provided in App. A) links the stability properties of the numerical schemes written in the form (27) with the spectral radius (i.e. the modulus of the largest eigenvalue, denoted by $\rho(\cdot)$) of the Jacobian matrix of the iteration function $\Phi$:

**Proposition 1.** *Let us consider the numerical scheme (27) and assume that $\Phi \colon \mathbb{R}^n \times \mathbb{R}^+ \times \mathbb{R}^+ \to \mathbb{R}^n$ is differentiable with respect to its first argument. If the condition*

$$\exists \, \Delta t_0 > 0, \, \exists \, \alpha \in \mathbb{R} : \quad \rho\left(\nabla_{\boldsymbol{\psi}} \Phi(\boldsymbol{\psi}, t, \Delta t)\right) \leq 1 + \alpha \, \Delta t \quad \forall \, \boldsymbol{\psi} \in \mathbb{R}^n, \, t \geq 0, \, \Delta t \in (0, \Delta t_0] \tag{29}$$

*holds true, then the scheme (27) is zero-stable. Moreover, if for a given $\Delta t$ the condition*

$$\exists \, \rho_0 < 1 : \quad \rho\left(\nabla_{\boldsymbol{\psi}} \Phi(\boldsymbol{\psi}, t, \Delta t)\right) \leq \rho_0 \quad \forall \, \boldsymbol{\psi} \in \mathbb{R}^n, \, t \geq 0, \tag{30}$$

*holds true, then the scheme is absolutely stable in correspondence to $\Delta t$.*

We remark that Prop. 1 can be used to study the stability of the monolithic, segregated and stabilized-segregated schemes for virtually any force generation model (not only for the minimal model (11)). However, focusing on the minimal model (11), we notice that in this case the considered schemes (i.e. (23), (24)–(25) and (24)–(26)) can be written in the following form, by writing $\boldsymbol{\psi}^{(k)} = (\mu_0^{(k)}, (\mathbf{y}^{(k)})^T)^T$, with

$\mathbf{y}^{(k)} = (\mu_1^{(k)}, \lambda^{(k)}, \lambda^{(k-1)})^T$:

$$\begin{cases} \mu_0^{(k+1)} = \dfrac{\mu_0^{(k)} + \mu_0^f \Delta t}{1 + r\,\Delta t} & k = 0, 1, \ldots \\ A\left(\mu_0^{(k+1)}, \Delta t\right) \mathbf{y}^{(k+1)} = B\left(\mu_0^{(k+1)}, \Delta t\right) \mathbf{y}^{(k)} + \mathbf{h}^{(k+1)} & k = 0, 1, \ldots \\ \mu_0^{(0)} = \mu_{0,0}, \qquad \mathbf{y}^{(0)} = \mathbf{y}_0, \end{cases} \quad (31)$$

where $\mathbf{h}^{(k+1)} = (\Delta t\, \mu_1^f, p(t_{k+1}), 0)^T$ and where the matrices $A$ and $B$ define the different schemes. Specifically, we have

$$A_{\mathrm{mon}}(\mu_0, \Delta t) = \begin{pmatrix} 1 + r\Delta t & -\mu_0 & 0 \\ a_{\mathrm{XB}} & K_p & 0 \\ 0 & 0 & 1 \end{pmatrix}, \quad B_{\mathrm{mon}}(\mu_0, \Delta t) = \begin{pmatrix} 1 & -\mu_0 & 0 \\ 0 & 0 & 0 \\ 0 & 1 & 0 \end{pmatrix}. \quad (32)$$

for the monolithic scheme (23). Then, for the segregated scheme (24)–(25), we have:

$$A_{\mathrm{segr}}(\mu_0, \Delta t) = \begin{pmatrix} 1 + r\Delta t & 0 & -\mu_0 \\ a_{\mathrm{XB}} & K_p & 0 \\ 0 & 0 & 1 \end{pmatrix}, \quad B_{\mathrm{segr}}(\mu_0, \Delta t) = \begin{pmatrix} 1 & 0 & -\mu_0 \\ 0 & 0 & 0 \\ 0 & 1 & 0 \end{pmatrix}. \quad (33)$$

Finally, the stabilized-segregated scheme (24)–(26) is obtained by modifying the matrices of the segregated scheme as follows:

$$\begin{aligned} A_{\text{stab-segr}}(\mu_0, \Delta t) &= A_{\mathrm{segr}}(\mu_0, \Delta t) + \begin{pmatrix} 0 & 0 & 0 \\ 0 & a_{\mathrm{XB}}\mu_0 & -a_{\mathrm{XB}}\mu_0 \\ 0 & 0 & 0 \end{pmatrix}, \\ B_{\text{stab-segr}}(\mu_0, \Delta t) &= B_{\mathrm{segr}}(\mu_0, \Delta t). \end{aligned} \quad (34)$$

Expression (27) is more general than (31), in the sense that the schemes that can be written in the form (31) represent a subset of those of (27). Hence, the notions of zero- and absolute stability (Defs. 1, 2) can be extended to the schemes written in the form (31), as a particular case of (27). Moreover, also Prop. 1 applies to the schemes written in the form (31), taking the following special form (see App. A for the proof):

**Proposition 2.** *Let us consider a numerical scheme written in the form (31). Let us suppose that, for any $\mu_0 \in [0, \mu_0^f/r]$ and for any $\Delta t > 0$*

*(H1) the matrix $A(\mu_0, \Delta t)$ is invertible;*

*(H2) the matrices $A(\mu_0, \Delta t)$ and $B(\mu_0, \Delta t)$ are differentiable with respect to $\mu_0$.*

*Let us denote $C(\mu_0, \Delta t) := A(\mu_0, \Delta t)^{-1} B(\mu_0, \Delta t)$. The scheme (31) is zero-stable, if*

$$\exists \Delta t_0 > 0, \exists \alpha \in \mathbb{R}: \quad \rho(C(\mu_0, \Delta t)) \leq 1 + \alpha\, \Delta t \qquad \forall \Delta t \in (0, \Delta t_0], \mu_0 \in [0, \mu_0^f/r]. \quad (35)$$

*Moreover, it is absolutely stable in correspondence to $\Delta t$ if*

$$\exists \rho_0 < 1: \qquad \rho(C(\mu_0, \Delta t)) < \rho_0 \qquad \forall \mu_0 \in [0, \mu_0^f/r]. \quad (36)$$

Conversely, it is easy to see that if $\rho(C(\mu_0, \Delta t)) > 1$ for some $\mu_0$ and $\Delta t$, the solution might blow up. Hence, the numerical stability of the method is determined by the spectral radius of the matrix $A^{-1}B$. Indeed, hypotheses (H1) and (H2) are easy to verify for each of the considered schemes. Let us then study the spectrum of this matrices $A^{-1}B$.

**Monolithic scheme.** The spectrum of the matrix $A_{\text{mon}}^{-1} B_{\text{mon}}$ contains the following eigenvalues:
$$\sigma_1 = \sigma_2 = 0, \quad \sigma_3 = \frac{K_p + a_{\text{XB}}\,\mu_0}{K_p + a_{\text{XB}}\,\mu_0 + \Delta t\, K_p\, r}.$$

Clearly, for physically meaningful values of the parameters (i.e. when they all take positive values), we have $\sigma_3 \in (0, 1)$. This entails that the monolithic scheme is zero-stable and unconditionally absolutely stable (i.e. without any restriction on $\Delta t$).

**Segregated scheme.** The spectrum of the matrix $A_{\text{segr}}^{-1} B_{\text{segr}}$ is given by:
$$\sigma_1 = 0, \quad \sigma_{2,3} = \sigma_\pm = \frac{K_p - a_{\text{XB}}\,\mu_0 \pm \sqrt{(K_p - a_{\text{XB}}\,\mu_0)^2 + 4 a_{\text{XB}}\,\mu_0\, K_p (1 + r\,\Delta t)}}{2 K_p (1 + r\,\Delta t)}.$$

It is easy to see that, for $\Delta t \to +\infty$ (this is simply a mathematical speculation, of course), we have $\sigma_\pm \to 0$, while, for $\Delta t \to 0$, we have $\sigma_+ \to 1$ and $\sigma_- \to -\frac{a_{\text{XB}}\,\mu_0}{K_p}$. Therefore, if $K_{\text{a}} = a_{\text{XB}}\,\mu_0 > K_p$ (i.e. when the stiffness of the active components is larger than the stiffness of the passive component), we have $\sigma_- < -1$ when $\Delta t$ is small, thus leading to nonphysical oscillations in the numerical solution. Therefore, whenever $K_{\text{a}} > K_p$, the segregated scheme cannot be absolutely stable when $\Delta t$ is lower than a given threshold. Moreover it is not even zero-stable, as for $\Delta t \to 0$ the spectral radius tends to a constant strictly lower than $-1$. Hence, since a consistent numerical scheme is convergent if and only if it is stable, the segregated scheme is never convergent when $K_{\text{a}} > K_p$. As a matter of fact, with the segregated scheme one cannot decrease the time step size $\Delta t$ aiming at achieving a prescribed accuracy in the numerical approximation, because this would compromise the stability of the method. This makes the segregated scheme somehow pointless, as a consistent numerical scheme is typically devised in such a way that, by progressively refining the discretization, the numerical solution converges to the exact one.

We notice that this result is in accordance with the results of [31], where it is shown that a necessary condition for stability of the segregated scheme is that the active stiffness should not be larger than the passive one. Similarly, in [34] the authors show that the segregated scheme may loose stability when the stiffness associated with the active stress is large, compared to the passive one, even if a sharp threshold is not derived.

**Stabilized-segregated scheme.** Finally, let us consider the stabilized-segregated scheme. The spectrum of the matrix $A_{\text{stab-segr}}^{-1} B_{\text{stab-segr}}$ contains the following eigenvalues:
$$\sigma_1 = \sigma_2 = 0, \quad \sigma_3 = \frac{K_p + a_{\text{XB}}\,\mu_0 + \Delta t\, a_{\text{XB}}\,\mu_0\, r}{K_p + a_{\text{XB}}\,\mu_0 + \Delta t\, a_{\text{XB}}\,\mu_0\, r + \Delta t\, K_p\, r}.$$

Clearly, we have $\sigma_3 \in (0, 1)$ for all $\Delta t > 0$. Therefore, the stabilized-segregated scheme, like the monolithic scheme, is zero-stable and absolutely unconditionally stable. In other words, the stabilization term has the effect of bringing the eigenvalue responsible for the lack of stability of the segregated scheme back into the unit interval, while preserving consistency.

In conclusion, we have the following result.

**Proposition 3.** *The monolithic scheme (23) and the stabilized-segregated scheme (24)–(26) are zero-stable and unconditionally absolutely stable.*

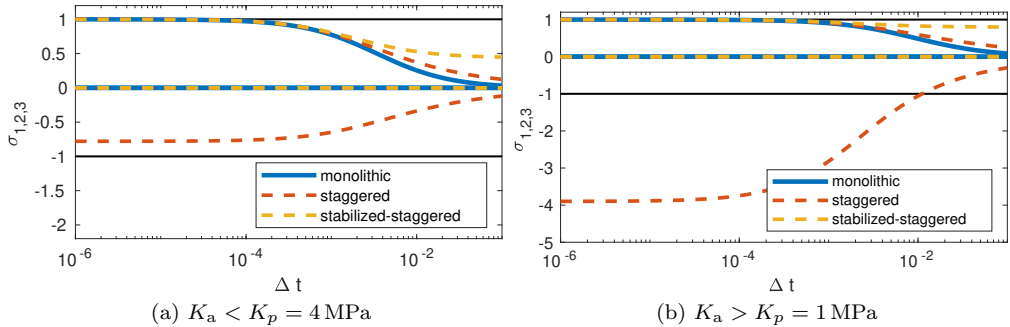

Figure 2: Quasistatic mechanics: eigenvalues associated with the three numerical schemes as a function of $\Delta t$, in the two cases $a_{\text{XB}} \mu_0 = K_a > K_p$ and $a_{\text{XB}} \mu_0 = K_a \leq K_p$.

In Fig. 2 we show the dependence on $\Delta t$ of the eigenvalues $\sigma_{1,2,3}$ associated with the three schemes, for realistic values of the parameters of the model (see App. B), where we employ the steady state value of the variable $\mu_0$, given by $\mu_0 = \mu_0^f/r$. In Fig. 2a we take $K_p = 4\,\text{MPa}$, so that $K_p > K_a$. In this case, the eigenvalues are always contained in $(-1, 1)$, which means that each of the three schemes is unconditionally absolutely stable. Conversely, in Fig. 2b, we take $K_p = 1\,\text{MPa}$, so that $K_p < K_a$. In this case, as we expected, the segregated scheme features an eigenvalue below $-1$ for $\Delta t$ not sufficiently small, while the others schemes are absolutely stable for any value of $\Delta t$.

In order to numerically validate the theoretical results derived above, we perform the following test. We employ the three schemes to approximate the solution of problem (22) in the quasistatic case (for simplicity we consider the case of no external load, by setting $p(t) \equiv 0$). In Fig. 3a we show the results obtained in the case $K_p = 4\,\text{MPa}$ (i.e. with a passive stiffness larger than the largest attained active stiffness). As expected, none of the three methods yield numerical oscillations. In Fig. 3b, instead, we consider $K_p = 1\,\text{MPa}$. In this case, after a while the active stiffness $K_a(t) = a_{\text{XB}} \mu_0(t)$ exceeds the passive one $K_p$, represented by a dashed black line. In the solution of the segregated scheme this leads, as expected, to nonphysical oscillations, that are successfully removed by the stabilization term introduced in the stabilized-segregated scheme.

### 4.3.2 Elastodynamics case

Consider now the case when the quasistatic approximation considered before is no longer assumed: inertia and viscosity will be restored in our mechanical model (see Eq. (22)). The previously considered numerical schemes can be easily generalized by approximating the additional terms with BDF formulas. More precisely, the term $K_p \lambda^{(k+1)}$ in the left-hand sides of Eqs. (23), (25) and (26) is replaced by

$$M \frac{\lambda^{(k+1)} - 2\lambda^{(k)} + \lambda^{(k-1)}}{\Delta t^2} + \sigma \frac{\lambda^{(k+1)} - \lambda^{(k)}}{\Delta t} + K_p \lambda^{(k+1)}. \tag{37}$$

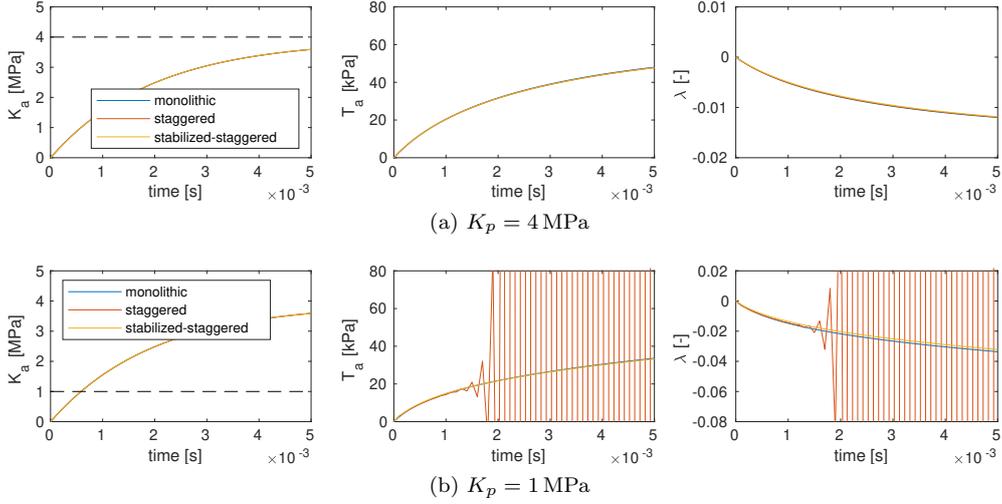

Figure 3: Quasistatic mechanics: results of simulations obtained with two different values of $K_p$ (represented by the dashed black line).

By introducing the following matrices

$$A_{\text{mass}} = B_{\text{mass}} = \frac{M}{\Delta t^2} \begin{pmatrix} 0 & 0 & 0 \\ 0 & 1 & -1 \\ 0 & 0 & 0 \end{pmatrix}, \qquad A_{\text{damp}} = B_{\text{damp}} = \frac{\sigma}{\Delta t} \begin{pmatrix} 0 & 0 & 0 \\ 0 & 1 & 0 \\ 0 & 0 & 0 \end{pmatrix},$$

the algebraic interpretation of the three schemes can be still written in the matrix form (31); however, the matrices reported in (32), (33) and (34) becomes, respectively:

$$A_{\text{mon}}^{\text{total}} = A_{\text{mon}} + A_{\text{mass}} + A_{\text{damp}}, \qquad B_{\text{mon}}^{\text{total}} = B_{\text{mon}} + B_{\text{mass}} + B_{\text{damp}},$$
$$A_{\text{segr}}^{\text{total}} = A_{\text{segr}} + A_{\text{mass}} + A_{\text{damp}}, \qquad B_{\text{segr}}^{\text{total}} = B_{\text{segr}} + B_{\text{mass}} + B_{\text{damp}},$$
$$A_{\text{stab-segr}}^{\text{total}} = A_{\text{stab-segr}} + A_{\text{mass}} + A_{\text{damp}}, \qquad B_{\text{stab-segr}}^{\text{total}} = B_{\text{stab-segr}} + B_{\text{mass}} + B_{\text{damp}}.$$

In Fig. 4 we report the eigenvalues (that corresponds to those displayed in Fig. 2b for the quasistatic case) when a non-null viscosity is added to the model (i.e. with $\sigma \neq 0$ and $M = 0$). The results show that the viscosity has a stabilizing effect on the segregated scheme for small values of $\Delta t$; however, there is still an interval of values of $\Delta t$ for which the segregated scheme is not absolutely stable. Conversely, both the monolithic and the stabilized-segregated schemes are unconditionally absolutely stable. Moreover, we notice that, for $\Delta t \to 0$ the eigenvalues associated with the stabilized-segregated scheme show a better adherence with those of the monolithic scheme than those of the segregated scheme.

Finally, we introduce also the inertia term (see Fig. 5). In this case, all the schemes feature a real eigenvalue and a pair of conjugate complex ones. Similarly to the previous case, the segregated scheme is absolutely unstable within an interval of values of $\Delta t$, while both the monolithic and the stabilized-segregated schemes are absolutely stable for any value of $\Delta t$. We remark that when only inertia is considered (i.e. when $\sigma = 0$ and $M \neq 0$), the results are qualitatively equivalent to those of Fig. (5).

To sum up, the following conclusions can be drawn.

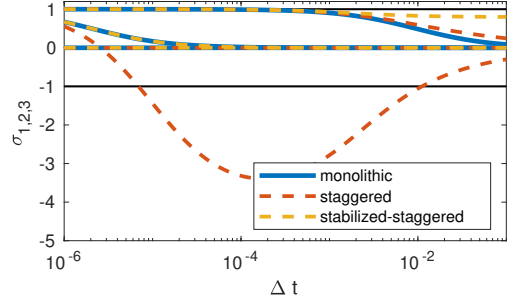

Figure 4: Non quasistatic mechanics: eigenvalues associated with the three numerical schemes as a function of $\Delta t$, in the case $\sigma = 10\,\mathrm{Pa\,s}$ and $M = 0$, with $a_{\mathrm{XB}}\,\mu_0 = K_{\mathrm{a}} > K_p$.

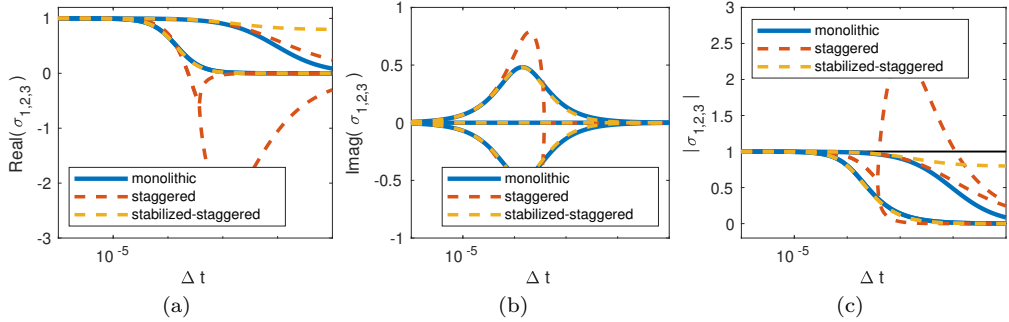

Figure 5: Non quasistatic mechanics: eigenvalues associated with the three numerical schemes as a function of $\Delta t$, in the case $\sigma = 10\,\mathrm{Pa\,s}$ and $M = 0.1\,\mathrm{Pa\,s^2}$, with $a_{\mathrm{XB}}\,\mu_0 = K_{\mathrm{a}} > K_p$. The three figures show respectively the real part, the imaginary part and the modulus of the eigenvalues.

- The monolithic scheme is always unconditionally absolutely stable.
- When $K_\mathrm{a} > K_p$, the segregated scheme may feature an eigenvalue outside the unit circle. In the quasistatic case, this happens for small values of $\Delta t$, whereas, when mass and/or viscosity are added to the mechanical model, it happens within an interval of values of $\Delta t$.
- Both in the quasistatic case and when mass and/or viscosity are considered, the effect of the stabilization term is that of bringing the eigenvalues responsible for the instability of the segregated scheme back into the unit circle, making the scheme absolutely stable for any value of $\Delta t$.

## 4.4 Interpretation as a fractional step scheme

In this section we show that the stabilized-segregated scheme can be interpreted as a fractional step scheme (see e.g. [38]), where the differential operator associated with the force generation model (11) is split into the sum of a term independent of the strain rate $\dot\lambda(t)$ and a velocity-dependent term; the latter is then associated with the operator of the mechanics model.

For simplicity, let us consider a constant $\mu_0$, thus focusing on the dynamics of $\mu_1$ and $\lambda$ only (i.e. on $(22)_2$–$(22)_3$). The stabilized-segregated scheme then reads

$$\begin{cases} \dfrac{\mu_1^{(k+1)} - \mu_1^{(k)}}{\Delta t} = \mu_1^f - r\,\mu_1^{(k+1)} + \dfrac{\lambda^{(k)} - \lambda^{(k-1)}}{\Delta t}\mu_0 & k = 0, 1, \ldots \\ M\dfrac{\lambda^{(k+1)} - 2\lambda^{(k)} + \lambda^{(k-1)}}{\Delta t^2} + \sigma\dfrac{\lambda^{(k+1)} - \lambda^{(k)}}{\Delta t} + K_p\lambda^{(k+1)} \\ \qquad = p(t_{k+1}) - a_\mathrm{XB}\left[\mu_1^{(k+1)} + \mu_0\left(\lambda^{(k+1)} - \lambda^{(k)}\right)\right] & k = 0, 1, \ldots \\ \widetilde\mu_1^{(0)} = \mu_{1,0}, \quad \lambda^{(-1)} = \lambda^{(0)} = \lambda_0 \end{cases} \tag{38}$$

We now write the differential equation for the variable $\mu_1$ as (see $(22)_2$)

$$\dot\mu_1(t) = \mu_1^f - r\,\mu_1(t) + \dot\lambda(t)\,\mu_0 =: \mathcal{L}(\mu_1(t), \dot\lambda(t)).$$

Then, we split the operator $\mathcal{L}$ as $\mathcal{L} = \mathcal{L}_1 + \mathcal{L}_2$, where

$$\mathcal{L}_1(\mu_1) = \mu_1^f - r\,\mu_1, \qquad \mathcal{L}_2(\dot\lambda) = \dot\lambda\,\mu_0.$$

Let us consider the following fractional-step scheme for the solution of $(22)_2$–$(22)_3$:

$$\begin{cases} \dfrac{\widetilde\mu_1^{(k+1,*)} - \widetilde\mu_1^{(k)}}{\Delta t} = \mathcal{L}_1(\widetilde\mu_1^{(k+1,*)}) = \mu_1^f - r\,\widetilde\mu_1^{(k+1,*)} & k = 0, 1, \ldots \\ \dfrac{\widetilde\mu_1^{(k+1)} - \widetilde\mu_1^{(k+1,*)}}{\Delta t} = \mathcal{L}_2\left(\dfrac{\lambda^{(k+1)} - \lambda^{(k)}}{\Delta t}\right) = \dfrac{\lambda^{(k+1)} - \lambda^{(k)}}{\Delta t}\mu_0 & k = 0, 1, \ldots \\ M\dfrac{\lambda^{(k+1)} - 2\lambda^{(k)} + \lambda^{(k-1)}}{\Delta t^2} + \sigma\dfrac{\lambda^{(k+1)} - \lambda^{(k)}}{\Delta t} + K_p\lambda^{(k+1)} \\ \qquad = p(t_{k+1}) - a_\mathrm{XB}\,\widetilde\mu_1^{(k+1)} & k = 0, 1, \ldots \\ \widetilde\mu_1^{(0)} = \mu_{1,0}, \quad \lambda^{(-1)} = \lambda^{(0)} = \lambda_0 \end{cases} \tag{39}$$

In this scheme we first advance the state $\widetilde\mu_1^{(k)}$ to an intermediate state $\widetilde\mu_1^{(k+1,*)}$, through the operator $\mathcal{L}_1$. Then, we advance from $\widetilde\mu_1^{(k+1,*)}$ to $\widetilde\mu_1^{(k+1)}$, through the operator $\mathcal{L}_2$.

The latter step is coupled in a monolithic manner with the mechanics equation. We have the following result.

**Proposition 4.** *The stabilized-segregated scheme* (38) *is equivalent to the fractional-step scheme* (39) *for the solution of* $(22)_2$–$(22)_3$. *More precisely, the two schemes produce the same sequence* $\lambda^{(k)}$, *for* $k \geq 0$, *while, for what concerns the discretization of the variable* $\mu_1$, *we have* $\mu_1^{(k)} = \widetilde{\mu}_1^{(k,*)}$, *for* $k \geq 1$.

*Proof.* The thesis follows from the equalities $\widetilde{\mu}_1^{(k)} = \mu_1^{(k)} + [\lambda^{(k)} - \lambda^{(k-1)}]\mu_0$ and $\widetilde{\mu}_1^{(k,*)} = \mu_1^{(k)}$. □

This equivalence provides a new interpretation of the stabilized-segregated scheme. Indeed, we notice that also the standard segregated scheme (24)–(25) can be interpreted as a fractional step method, where the first step is associated with the full operator $\mathcal{L}$ and the second step is associated with the operator of the mechanical problem only. Hence, the stabilized-segregated scheme can be interpreted as a modification of the standard staggered scheme, where the velocity-dependent part of the operator $\mathcal{L}$ – associated with the velocity-related feedback loop and thus responsible for the lack of stability of the segregated scheme – is solved simultaneously with the mechanics model, thus eliminating such instabilities. We remark that this is just an interpretation of the stabilized-segregated scheme, in light of the equivalence shown above. In practice, one does not have to actually couple the two models in a monolithic way: the stabilized-segregated scheme preserves its segregated feature.

# 5 Stabilized-segregated scheme for a generic force generation model

In the previous sections, we have derived, based on energetic considerations, a fully segregated but stable numerical scheme to couple active mechanics models with the distribution-moments force generation model of Eq. (11). Despite its derivation is rooted into the physics of this specific model, this scheme can be easily generalized to other force generation models.

## 5.1 Generalization to different models

In our stabilized-segregated scheme, the extra term $K_\mathrm{a}^{(k+1)}\left(\lambda^{(k+1)} - \lambda^{(k)}\right)$ was added to $T_\mathrm{a}^{(k+1)}$ in the mechanical model. This suggests the way to generalize this scheme to other force generation models, provided we suitably define a value for the active stiffness $K_\mathrm{a}$.

A formula for $K_\mathrm{a}$ can be based upon physical considerations, similarly to what done in Sec. 3.2. Alternatively, we can formally proceed as follows. In abstract terms, the active stiffness is defined as

$$K_\mathrm{a} = \frac{\partial \dot{T}_\mathrm{a}}{\partial \dot{\lambda}}. \tag{40}$$

Let us start from the generic expression of Eq. (1). Then,

$$\dot{T}_\mathrm{a} = \nabla_\mathbf{r} g \cdot \mathbf{h}, \tag{41}$$

which entails

$$K_\mathrm{a} = \frac{\partial \dot{T}_\mathrm{a}}{\partial \dot{\lambda}} = \nabla_\mathbf{r} g \cdot \frac{\partial \mathbf{h}}{\partial \dot{\lambda}}. \tag{42}$$

| Abbreviation | # variables | Family (velocity-dependence) | Ref. |
|---|---|---|---|
| NHS06 | 5 | Fading-memory models | [30] |
| L17 | 6 | Distortion-decay models | [27] |
| RDQ20-MF | 20 | Physics-based models (Huxley formalism) | [43] |

Table 1: List of the the force generation models used in this paper.

In conclusion, we can use Eq. (42) to obtain the active stiffness corresponding to any activation model written in the form of Eq. (1).

*Remark* 1. Not all the models available in the literature fall in the class of Eq. (1). For some of them, for instance, the right-hand side $\mathbf{h}$ depends on $T_\mathrm{a}$ itself. However, this modification does not affect the definition of Eq. (42). In some other models, instead, the active tension explicitly depends on $\lambda$ in the following manner:

$$T_\mathrm{a}(t) = g(\mathbf{r}(t), \lambda(t)). \tag{43}$$

In this case, by proceeding formally, one gets:

$$K_\mathrm{a} = \frac{\partial \dot{T}_\mathrm{a}}{\partial \dot{\lambda}} = \nabla_\mathbf{r} g \cdot \frac{\partial \mathbf{h}}{\partial \dot{\lambda}} + \frac{\partial g}{\partial \lambda}. \tag{44}$$

For all the models considered in this paper, we have tested both options (Eq. (42) and Eq. (44)), without obtaining significantly different results. As a matter of fact, typically the additional term $\frac{\partial g}{\partial \lambda}$ is negligible compared to the other term. For this reason, Eq. (42) should be preferred to Eq. (44) because of its simplest implementation.

In the next sections we apply this procedure to some of the force generation models available in the literature and currently used in multiscale cardiac electromechanical simulations. We consider three among the most popular families of models, classified according to the formalism used to describe the velocity-dependent feedback on the force generating machinery (see Tab. 1). Specifically, we consider two phenomenological families (fading-memory models and distortion-decay models) and the class of physics-based models. For each family we consider a representative model, and we provide the expression of its active stiffness $K_\mathrm{a}$. Then, in Sec. 5.2, we present some numerical results obtained with these models.

### 5.1.1 Fading-memory models

Fading-memory models [3, 20, 21, 28, 30] are based on a phenomenological description of the velocity-related effects on the generated force. In these models, a nonlinear function of the ratio between the active tension and the isometric tension $T_\mathrm{a}^\mathrm{iso}$ is set equal to an integral convolution of the past history of the strain rate $\dot{\lambda}$ with a nonlinear kernel $\varphi$:

$$Q(T_\mathrm{a}/T_\mathrm{a}^\mathrm{iso}) = \int_{-\infty}^{t} \varphi(t-\tau)\dot{\lambda}(\tau)d\tau. \tag{45}$$

The nonlinear kernel $\varphi$ is set equal to a finite sum of negative exponentials:

$$\varphi(t) = \sum_{i=1}^{N} A_i e^{-\alpha_i t},$$

so that the current tension is more influenced by the most recent length changes than by the earlier ones. Then, the function $Q$ is defined so that the constant-velocity solution fits the Hill curve [18].

Let us consider the Niederer-Hunter-Smith model [30] (henceforth denoted as NHS06 model), where the velocity-dependent effects are accounted for by means of the fading-memory formalism. This model is described in detail in App. (C). In this discussion, it is sufficient to know that the force predicted by this model is given by

$$T_\mathrm{a} = T_\mathrm{ref} \left(1 + \beta_0\, \lambda\right) \frac{z}{z_\mathrm{max}(\lambda)}\, K\left(\sum_{i=1}^{3} Q_i\right),$$

where $T_\mathrm{ref}$ is a reference tension, $\beta_0$ is a coefficient associated with the length-dependence of the steady-state tension, $z/z_\mathrm{max}(\lambda)$ represents the fraction of available binding sites (with respect to the maximum binding sites recruitable for a given sarcomere length) and, finally, $K$ is the velocity-dependent function, obtained by inverting the function $Q$ of Eq. (45), whose definition if provided in App. (C). The evolution of the phenomenological state variables $Q_i$ is described by the following ODE:

$$\dot{Q}_i = A_i\, \dot{\lambda} - \alpha_i Q_i, \qquad i = 1, 2, 3. \tag{46}$$

Hence, by formally applying Eq. (42) to the NHS06 model, we obtain the following definition for the active stiffness:

$$K_\mathrm{a} = T_\mathrm{ref} \left(1 + \beta_0\, \lambda\right) \frac{z}{z_\mathrm{max}(\lambda)}\, K'\left(\sum_{i=1}^{3} Q_i\right) \sum_{i=1}^{3} A_i, \tag{47}$$

where $K'$ denotes the derivative of $K$ (its definition is provided in App. (C)).

### 5.1.2 Distortion-decay models

In distortion-decay models [13, 40, 46], the population of attached XBs is split into a number of groups. A variable describing the average distortion of the XBs is associated with each group, so that the generated force can be computed as the product of the XB stiffness and the XB distortion variables, respectively weighted by the size of each group.

For instance, in the model proposed by Land and coworkers in [27] (see App. D), henceforth denoted as L17 model, XBs are split into weakly-bound (i.e. in the pre-powerstroke configuration) and strongly-bound (i.e. in the post-powerstroke configuration). The average distortions associated with these groups are respectively tracked by the variables $\zeta_w(t)$ and $\zeta_s(t)$, which evolve according to the following phenomenological model:

$$\begin{cases} \dot{\zeta}_w = A_w \dot{\lambda} - c_w\, \zeta_w, \\ \dot{\zeta}_s = A_s \dot{\lambda} - c_s\, \zeta_s. \end{cases} \tag{48}$$

Finally, the active tension is computed as

$$T_\mathrm{a} = h(\lambda) \frac{T_\mathrm{ref}}{r_s} \left[(1 + \zeta_s)\, S + \zeta_w\, W\right], \tag{49}$$

where $h(\lambda)$ models the force-length relationship, $T_\mathrm{ref}$ is a reference tension, $r_s$ is the steady-state duty-ratio, while $S$ and $W$ are state variables representing the fraction of XBs in weakly-bound and strongly-bound configuration, respectively.

From Eq. (42), it follows that the active stiffness for the L17 model reads

$$K_\mathrm{a} = h(\lambda) \frac{T_\mathrm{ref}}{r_s} \left[A_s\, S + A_w\, W\right]. \tag{50}$$

### 5.1.3 Physics-based models

When we consider physics-based models of active force generations, the active stiffness $K_\mathrm{a}$ may be derived by means of physics-driven considerations, similarly to what done in Sec. 3.2 for the distribution-moments model. Alternatively, a formal derivation based on Eq. (42) can be employed as well.

The mean-field model that we proposed in [43] (see App. E), henceforth denoted as RDQ20-MF model, is a physics-based model where the calcium-driven activation is described within the formalism of [41], simplified by a mean-field approximation, while XB cycling is described within the distribution-moments formalism. Specifically, the population of binding sites is split into two groups, according to the state of the associated tropomyosin (either permissive or not). The active tension predicted by this model is hence given by

$$T_\mathrm{a} = a_\mathrm{XB}\, \chi_\mathrm{so}(\lambda) \left[\mu_\mathcal{P}^1 + \mu_\mathcal{N}^1\right], \tag{51}$$

where $\chi_\mathrm{so}(\lambda)$ is the length-dependent single-overlap ratio, while $\mu_\mathcal{P}^p$ and $\mu_\mathcal{N}^p$ are the $p$-th order distribution-moments associated with the permissive and non-permissive groups, respectively. By formally applying Eq. (42), we get that the active stiffness is given by

$$K_\mathrm{a} = a_\mathrm{XB}\, \chi_\mathrm{so}(\lambda) \left[\mu_\mathcal{P}^0 + \mu_\mathcal{N}^0\right]. \tag{52}$$

A thorough physics-based derivation, yielding the same result, can be found in [45].

Similarly, the active stiffness of models belonging to the H57 family in its original form (i.e. without the distribution-moments reduction) can be computed as

$$K_\mathrm{a} = \frac{a_\mathrm{XB}}{D_M} \int_{-\infty}^{+\infty} n(x,t)dx. \tag{53}$$

Another family of physics-based models is that of Markov Chain models [22, 53, 54]. In these models, a large number of realizations of a stochastic process describing the dynamics of regulatory and contractile proteins is simulated in parallel, and the obtained results are averaged. For these models the active stiffness can be obtained on a physical basis, computing the total stiffness of attached XBs as the average of the realizations of the stochastic processes. This is equivalent to computing the integral of Eq. (53) by Monte Carlo sampling.

## 5.2 Numerical results

In this section, we show the results obtained by applying the proposed stabilized numerical scheme to the three models considered in Sec. 5.1 and listed in Tab. 1, with the goal of showing its effectiveness.

Similarly to Sec. 4.3, we consider a zero-dimensional model for the tissue strain in the myofibers direction. Specifically, we consider the following model:

$$M\ddot{\lambda}(t) + \sigma\dot{\lambda}(t) + \frac{\partial \mathcal{W}}{\partial \lambda}(\lambda(t)) + T_\mathrm{a}(t) = p(t), \quad t \geq 0, \tag{54}$$

with $\sigma = 10\,\mathrm{Pa\,s}$ and $M = 0.1\,\mathrm{Pa\,s}^2$ and where $\mathcal{W}(\lambda)$ denotes an elastic potential. In place of the quadratic potential $\mathcal{W}(\lambda) = \frac{1}{2} K_p \lambda^2$ considered in Sec. 4.3, in this section we consider the nonlinear elastic potential $\mathcal{W}(\lambda) = \frac{1}{2} K_p \lambda \log(1 + \lambda)$, with $K_p = 1\,\mathrm{MPa}$.

Hence, the coupled model of active mechanics reads as follows:

$$\begin{cases} \dot{\mathbf{r}}(t) = \mathbf{h}\left(\mathbf{r}(t), [\mathrm{Ca}^{2+}]_\mathrm{i}(t), \lambda(t), \dot{\lambda}(t)\right), & t \geq 0, \\ T_\mathrm{a}(t) = g\left(\mathbf{r}(t), \lambda(t)\right), & t \geq 0, \\ K_\mathrm{a}(t) = q\left(\mathbf{r}(t), \lambda(t)\right), & t \geq 0, \\ M \ddot{\lambda}(t) + \sigma \dot{\lambda}(t) + \dfrac{\partial \mathcal{W}}{\partial \lambda}(\lambda(t)) = p(t) - T_\mathrm{a}(t), & t \geq 0, \\ \mathbf{r}(0) = \mathbf{r}_0, \quad \lambda(0) = \lambda_0, \quad \dot{\lambda}(0) = 0, \end{cases} \qquad (55)$$

where, for the sake of generality, we consider the case when the active tension explicitly depends on $\lambda(t)$. We consider the following numerical scheme:

$$\begin{cases} \dfrac{\mathbf{r}^{(k+1)} - \mathbf{r}^{(k)}}{\Delta t} = \widehat{\mathbf{h}}\left(\mathbf{r}^{(k)}, \mathbf{r}^{(k+1)}, [\mathrm{Ca}^{2+}]_\mathrm{i}^{(k+1)}, \lambda^{(*)}, \dot{\lambda}^{(*)}\right), & k \geq 0, \\ T_\mathrm{a}^{(k+1)} = g\left(\mathbf{r}^{(k+1)}, \lambda^{(*)}\right), & k \geq 0, \\ K_\mathrm{a}^{(k+1)} = q\left(\mathbf{r}^{(k+1)}, \lambda^{(*)}\right), & k \geq 0, \\ M \dfrac{\lambda^{(k+1)} - 2\lambda^{(k)} + \lambda^{(k-1)}}{\Delta t^2} + \sigma \dfrac{\lambda^{(k+1)} - \lambda^{(k)}}{\Delta t} + \dfrac{\partial \mathcal{W}}{\partial \lambda}(\lambda^{(k+1)}) = p^{(k+1)} - T_\mathrm{a}^{(*)}, & k \geq 0, \\ \mathbf{r}^{(0)} = \mathbf{r}_0, \quad \lambda^{(-1)} = \lambda^{(0)} = \lambda_0 \end{cases}$$
(56)

The numerical scheme used to treat the core block of force generation is defined through the function $\widehat{\mathbf{h}}$, that must satisfy, for consistency: $\widehat{\mathbf{h}}(\mathbf{r}, \mathbf{r}, c, \lambda, v) = \mathbf{h}(\mathbf{r}, c, \lambda, v)$. In particular, for $\widehat{\mathbf{h}}(\mathbf{r}_1, \mathbf{r}_2, c, \lambda, v) = \mathbf{h}(\mathbf{r}_1, c, \lambda, v)$ we have a fully-explicit scheme while for $\widehat{\mathbf{h}}(\mathbf{r}_1, \mathbf{r}_2, c, \lambda, v) = \mathbf{h}(\mathbf{r}_2, c, \lambda, v)$ we have a fully-implicit scheme. More in general, the function $\widehat{\mathbf{h}}$ defines a semiimplicit (or implicit-explicit, IMEX) scheme. Semiimplicit schemes allow to explicitly treat some nonlinearities, making the scheme linear, while keeping implicit those dependencies that would otherwise lead to instabilities. The numerical schemes used to treat the core models of force generation are described in Apps. C, D and E.

Concerning the coupling between the models of force generation and of tissue mechanics, the three schemes considered in what follows are:

- monolithic scheme, defined as:

$$\lambda^{(*)} = \lambda^{(k+1)}, \quad \dot{\lambda}^{(*)} = \frac{\lambda^{(k+1)} - \lambda^{(k)}}{\Delta t}, \quad T_\mathrm{a}^{(*)} = T_\mathrm{a}^{(k+1)};$$

- segregated scheme, defined as:

$$\lambda^{(*)} = \lambda^{(k)}, \quad \dot{\lambda}^{(*)} = \frac{\lambda^{(k)} - \lambda^{(k-1)}}{\Delta t}, \quad T_\mathrm{a}^{(*)} = T_\mathrm{a}^{(k+1)};$$

- stabilized-segregated scheme, defined as:

$$\lambda^{(*)} = \lambda^{(k)}, \quad \dot{\lambda}^{(*)} = \frac{\lambda^{(k)} - \lambda^{(k-1)}}{\Delta t}, \quad T_\mathrm{a}^{(*)} = T_\mathrm{a}^{(k+1)} + K_\mathrm{a}^{(k+1)} \left(\lambda^{(k+1)} - \lambda^{(k)}\right).$$

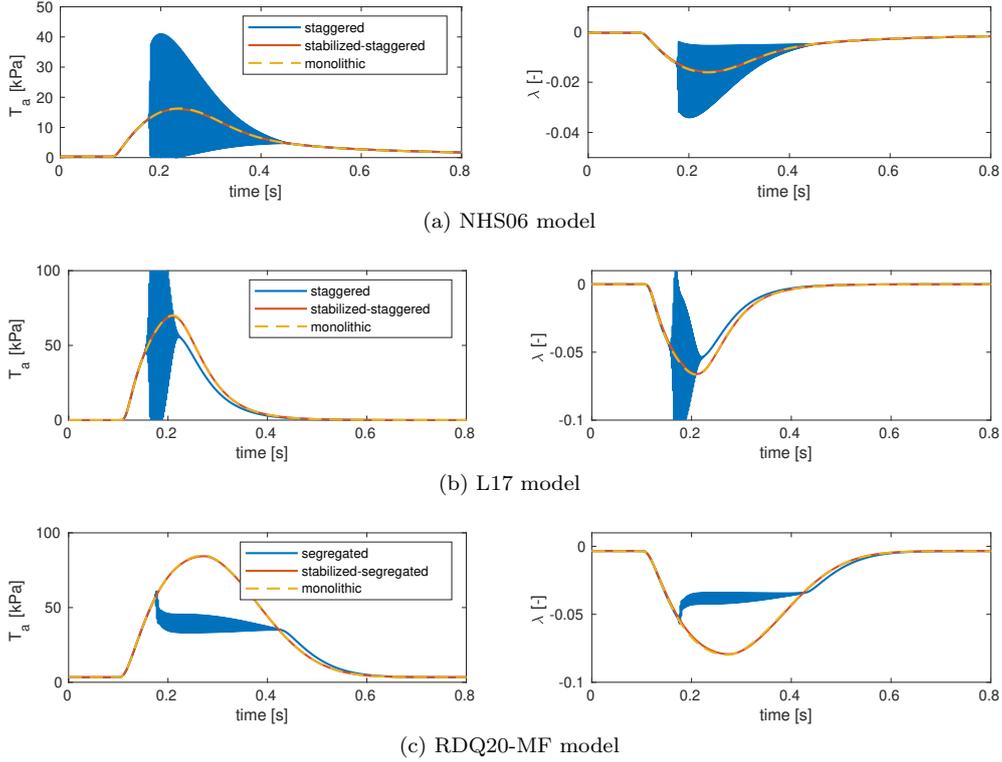

Figure 6: *Test Case 1*: Comparison of the results ($T_\mathrm{a}$, left and $\lambda$, right, versus time) obtained with the monolithic, segregated and stabilized-segregated schemes for three different force generation models.

### 5.2.1 *Test Case 1*: isotonic twitch

The first test case that we consider is a fiber contracting under isotonic conditions, that is against a constant load $p(t) \equiv 0\,\mathrm{Pa}$. To trigger the contraction of the fiber, we consider the following idealized calcium transient [52]:

$$[\mathrm{Ca}^{2+}]_\mathrm{i}(t) = c_0 + \frac{c_\mathrm{max} - c_0}{\beta} \left[ e^{-\frac{t-t_0}{\tau_1}} - e^{-\frac{t-t_0}{\tau_2}} \right] \mathbb{1}_{t \geq t_0}, \tag{57}$$

where

$$\beta = \left(\frac{\tau_1}{\tau_2}\right)^{-\left(\frac{\tau_1}{\tau_2}-1\right)^{-1}} - \left(\frac{\tau_1}{\tau_2}\right)^{-\left(1-\frac{\tau_2}{\tau_1}\right)^{-1}},$$

with $c_\mathrm{max} = 1.6\,\mathrm{\mu M}$, $c_0 = 0.1\,\mathrm{\mu M}$, $t_0 = 0.1\,\mathrm{s}$, $\tau_1 = 0.02\,\mathrm{s}$, $\tau_2 = 0.05\,\mathrm{s}$.

In Fig. 6 we compare the results obtained for each model with the three numerical schemes. While the monolithic scheme is stable for each model, the segregated scheme yields to numerical oscillations during the phase of maximal activation. This is coherent with our theoretical analysis of Sec. 4.3, that shows that such numerical instabilities arise when the active stiffness overcomes the passive one, i.e. when the muscle tissue is activated beyond a given threshold. Finally, the stabilized-segregated scheme successfully removes these numerical artifacts, for each of the three models.

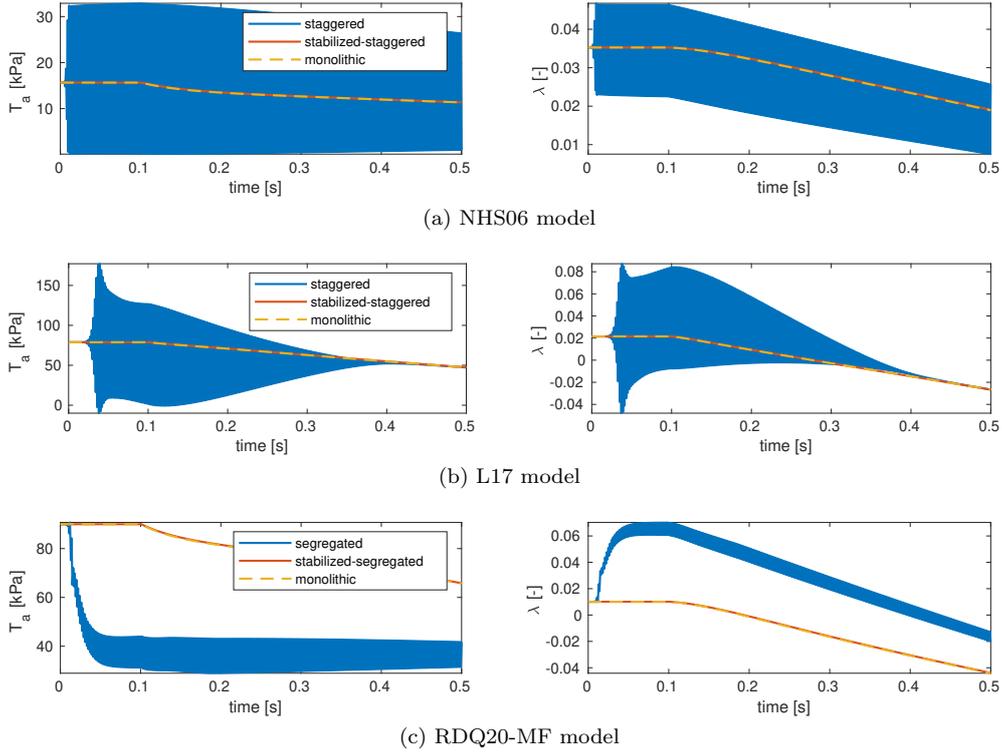

Figure 7: *Test Case 2*: Comparison of the results obtained with the monolithic, segregated and stabilized-segregated schemes for three different force generation models.

#### 5.2.2 *Test Case 2*: pressure ramp

In this second test case, we start from an equilibrium configuration where the muscle fiber is activated with a constant calcium concentration, against a constant external load $\bar{p}$. Then, we progressively decrease the external load, to simulate what happens during the ejection phase of a cardiac ventricle, yielding the the shortening of the muscle fiber.

Specifically, we consider a constant calcium concentration of $[\text{Ca}^{2+}]_i = 0.6\,\mu\text{M}$ for the NHS06 and L17 model and of $[\text{Ca}^{2+}]_i = 0.3\,\mu\text{M}$ for the RDQ20-MF model. We consider an initial load $\bar{p} = 50\,\text{kPa}$ for the NHS06 model and $\bar{p} = 100\,\text{kPa}$ for the L17 and RDQ20-MF models. Starting at time $t = 0.1\,\text{s}$, we linearly decrease $p(t)$ to reach zero in $0.5\,\text{s}$.

As expected, the numerical results (see Fig. 7) show that the shortening of the tissue leads to a decrease of active tension (due to the above mentioned force-velocity relationship). Similarly to *Test Case 1*, the stabilized-segregated scheme is successful in eliminating the nonphysical oscillations associated with the non-stabilized segregated scheme, thus giving a result consistent with that of the monolithic scheme.

#### 5.2.3 Convergence analysis

In order to verify the numerical consistency of the proposed stabilized-segregated scheme and to quantify the introduced error, we carry out a convergence analysis with respect to the time step length $\Delta t$. We consider again *Test Case 1* and, for each of the three force generation models of Tab. 1, we compute the numerical solution obtained with the monolithic and stabilized-segregated schemes for different values of $\Delta t$. Then, we compare the obtained time-evolution of the displacement with a reference solution. The latter, called $\lambda_{\text{ref}}$, in fact is obtained with the monolithic scheme for a very fine time step ($\Delta t = 10^{-4}$ s). To quantify the numerical error, we compute the $L^\infty$ and the mean $L^2$ distance from the reference solution $\lambda_{\text{ref}}$, interpolated on the finest time discretization:

$$
\begin{aligned}
e_2 &= \sqrt{\frac{1}{N_t} \sum_{k=0,\dots,N_t} \left|\lambda^{(k)} - \lambda^{(k)}_{\text{ref}}\right|^2}, \\
e_\infty &= \max_{k=0,\dots,N_t} \left|\lambda^{(k)} - \lambda^{(k)}_{\text{ref}}\right|,
\end{aligned}
\qquad (58)
$$

where $N_t$ denotes the number of time steps. In Fig. 8 we plot the errors for the three models, where we use, for the time-discretization of the force generation core models, either a semimplicit scheme, or an implicit scheme, or both of them.

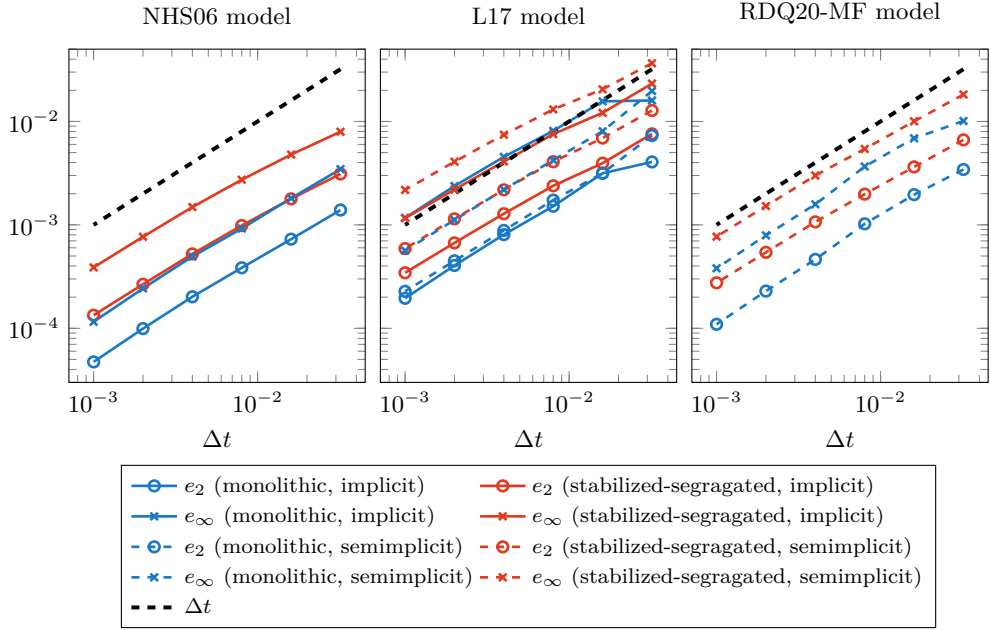

Figure 8: Convergence analysis with respect to the time step length $\Delta t$ for the three force generation models of Tab. 1.

First of all, we notice that, in all the considered cases, the numerical solutions obtained with the stabilized-segregated scheme converge to the exact solution as $\Delta t \to 0$. Indeed, the stabilized scheme that we propose introduces a perturbation of the segregated scheme by an infinitesimal term of order $\mathcal{O}(\Delta t)$ for $\Delta t \to 0$ (i.e. the scheme remains numerically consistent).

In all the cases, we numerically obtain, for the solution of the stabilized-segregated scheme, an order of convergence equal to one. Hence, the introduction of the stabilization term does not affect the overall first-order convergence of the coupled scheme of Eq. (56), due to the first-order finite difference schemes used to approximate time derivatives. In practice, for a given $\Delta t$, the errors obtained with the stabilized-segregated scheme are only slightly larger than the errors obtained with a monolithic scheme, for the same $\Delta t$.

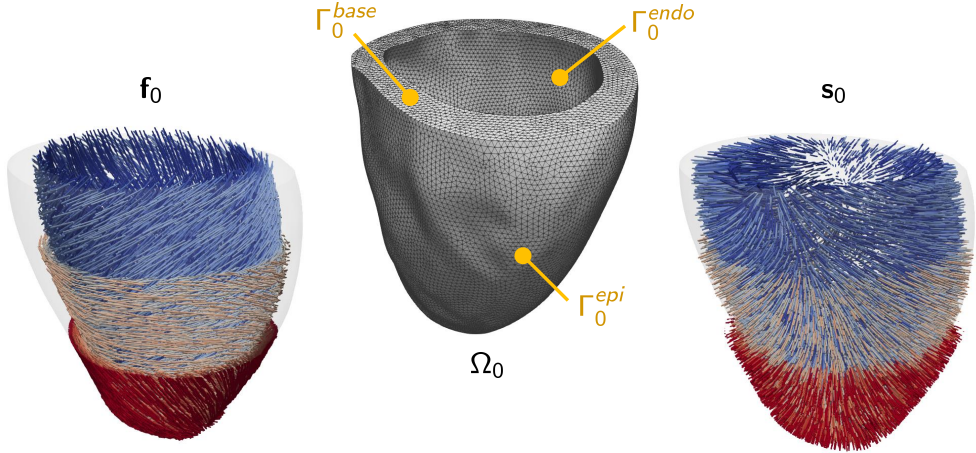

Figure 9: Computational mesh used for the three-dimensional electromechanical simulation of Sec. 5.2.4 (center) and reconstructed fibers ($\mathbf{f}_0$, left) and sheets ($\mathbf{s}_0$, right) fields.

We remark that in case higher order schemes are used for the time discretization of the core models, the proposed stabilization term should be generalized to be of the same order as the time stepping scheme, in order to preserve the overall order of convergence. However, one should check that this high-order generalization shares the stability properties of the first-order stabilized scheme proposed in this paper. This will be the object of a future work.

### 5.2.4 Three-dimensional left ventricle electromechanics

In this section we show the numerical results obtained by applying the proposed stabilized-segregated scheme in the context of the electromechanical modeling of a three-dimensional left ventricle. Specifically, we consider a computational domain derived from the Zygote heart model [59], which represents that of an average 21 years old healthy man. We generate the fibers and sheets direction within this domain by employing the rule-based algorithm proposed in [2]. The computational mesh and the fiber and sheets normal directions ($\mathbf{f}_0$ and $\mathbf{s}_0$, respectively) are represented in Fig. 9.

To model cardiac electrophysiology, we employ the monodomain equation [9, 10] and the ten Tusscher-Panfilov ionic model [48]. We describe the force generation phenomenon through the RDQ20-MF model (see App. E). We model the passive behaviour of the myocardium by means of the quasi-incompressible exponential constitutive law of [51]. At the epicardium ($\Gamma_0^{epi}$, see Fig. 9) we impose the generalized Robin boundary conditions of [15, 35], modeling the visco-eastic interaction of the cardiac walls with the pericardium. At the ventricular base ($\Gamma_0^{base}$, see Fig. 9), we set the energy-consistent boundary condition that we proposed in [44], which accounts for the effect of the cardiac muscle beyond the artificial boundary of the ventricular base. Finally, we couple the three-dimensional electromechanical ventricular model with a zero-dimensional Windkessel model of blood circulation [55]. The ventricular blood pressure of the zero-dimensional model provides the boundary condition at the endocardium to the three-dimensional model ($\Gamma_0^{endo}$, see Fig. 9). Further details on this electromechanical model can be found in [45].

Concerning the numerical solution of this model, we consider a tetrahedral computational mesh with $354 \cdot 10^3$ cells and $65 \cdot 10^3$ degrees of freedom. We use piece-wise linear Finite Elements to discretize in space the variables of the electromechanical model, while we employ first-order finite differences to discretize time derivatives, within the IMEX scheme described in [45]. The different core models are coupled following a fully segregated strategy [11], by orderly updating the solution of the monodomain equation, of the ionic model, of the force generation model and finally of the mechanical model. Moreover, to highlight the versatility of the proposed stabilized scheme, we consider a time-splitting scheme, by updating the solution of the mechanical subproblem every two iterations of the remaining physics. In other words, we employ a time step $\Delta t = 4 \cdot 10^{-4}$ s for the mechanical model and a time step $\Delta t = 2 \cdot 10^{-4}$ s for the other models.

The red lines in Fig. 10 show the results obtained with this segregated scheme, without stabilization. We observe the presence of spurious oscillations, due to the feedback of the tissue strain rate on the force generation dynamics. Moreover, we remark that, for larger time step lengths, the oscillations quickly lead to failure of the Newton solver. Conversely, the blue lines show the results obtained by adding the stabilization term proposed in this paper (see Eq. 19). These results show that the proposed scheme successfully achieves the stabilization of the observed spurious oscillations.

# 6 Conclusions

In this paper we proposed a new numerical scheme to couple force generation models with models describing tissue mechanics in cardiac electromechanics, aimed at curing numerical instabilities arising from the feedback that the tissue strain rate has on the force generation dynamics. The proposed numerical scheme is segregated, in the sense that the two subproblems are solved sequentially at each time step. This is more computationally attractive than a fully coupled (monolithic) scheme. Moreover, and most importantly, it is numerically stable, as we analytically proved for a model problem and we showed by means of several numerical tests. With respect to standard segregated coupling schemes, this scheme features an additional term, that can be seen as a numerically consistent stabilization term (of the order of $\mathcal{O}(\Delta t)$).

The proposed scheme is grounded on energetical considerations on the generation of the active stress at the microscale. Indeed, we showed that in the formula expressing the active stress generally employed in the literature (i.e. Eq. (7), see e.g. [16, 19]), a dependence on the tissue strain $\lambda$ is somehow hidden. This is due to a hybrid Lagrangian-Eulerian formalism by which the macroscopic displacement is typically expressed in the reference configuration, whereas the XB distortion is typically expressed with respect to the current frame of reference. Conversely, by framing all the variables in a more coherent fully Lagrangian reference system, we derived an expression for the active stress tensor, equivalent to the original one, which however leads to a different scheme when discretized with respect to time (see Eq. (19)). In the proposed scheme, the active tension is not expressed as a constant value, rather it depends on the tissue strain in the fibers direction, more coherently with the microscopic phenomenon that generates the macroscopic force.

To gain insight into the source of instabilities and to better appreciate the effects of the proposed stabilization term, we considered a minimal model of force generation and of mechanics. Then, we analyzed the numerical stability of different numerical schemes (monolithic, segregated and stabilized-segregated) by studying the eigenvalues of the corresponding transition matrices. We thus showed that the monolithic scheme is always absolutely stable, regardless of the time step size $\Delta t$ (that is, it is unconditionally absolutely stable). Conversely, when the active stiffness is larger than the passive stiffness in the fibers direction ($K_a > K_p$), the standard segregated scheme features eigenvalues outside the unit circle. In particular, in the case of quasistatic mechanics, the transition matrix associated with the segregated scheme has a real eigenvalue strictly lower than $-1$ (thus leading to nonphysical oscillations) for any $\Delta t$ lower than a critical value. This entails that this kind of instability is rather tough, as it cannot be removed by simply refining the temporal discretization. In other terms, the segregated scheme cannot be convergent when the active stiffness is larger than the passive stiffness in the fibers direction (as it happens during systole for realistic values of the parameters). Adding inertia and/or damping to the mechanical model improves the stability of the segregated scheme for the smallest values of $\Delta t$. However, for physiological values of the parameters, a wide window of values of $\Delta t$ for which the scheme cannot be absolutely stable is still present. We showed that the stabilization term proposed in this paper brings all the eigenvalues associated with the segregated scheme back into the unit circle, thus ensuring numerical stability. In conclusion, the proposed scheme is unconditionally absolutely stable, as the monolithic one.

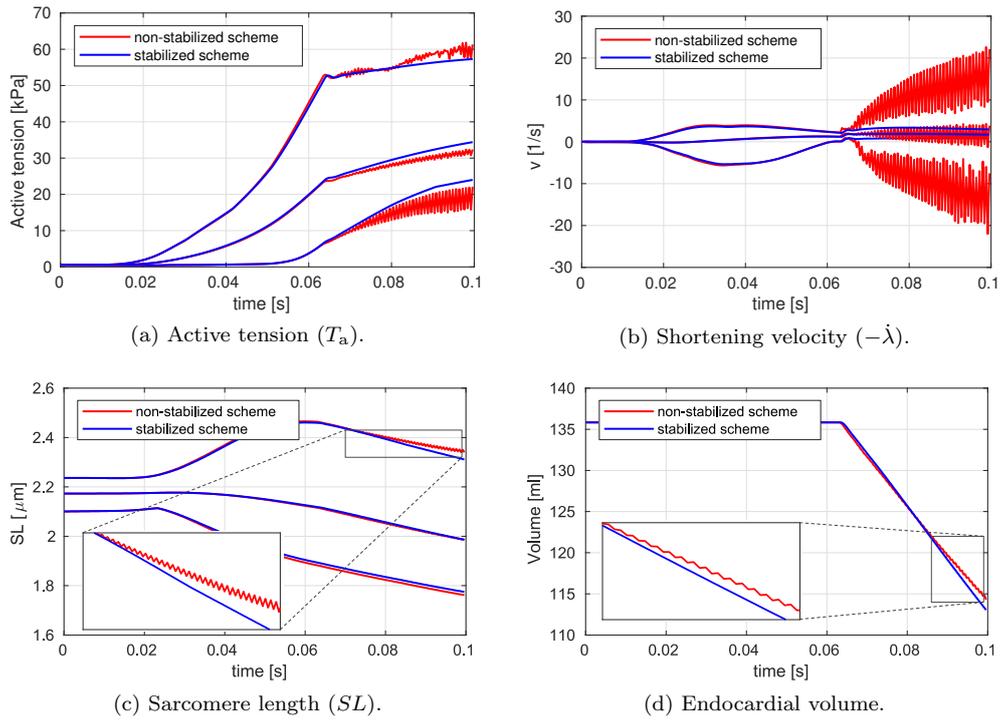

Figure 10: Results of multiscale cardiac electromechanis simulations obtained with the segregated scheme (red lines) and the stabilized-segregated scheme (blue lines). In (a)-(b)-(c) the three lines refer to the minimum, average and maximum value over the computational domain.

The proposed scheme can be easily generalized to several force generation models. Indeed, in order to apply this scheme to a given model, one only needs to derive a formula for the active stiffness associated with the model. This can be done either by formally applying a formula based on the model equations (Eq. (42)), or on the basis of a physics-driven analysis. As a matter of fact, we showed how the active stiffness can be derived for several families of force generation models available in the literature, including fading-memory models, distortion-decay models, and physics based models (e.g. models based on the Huxley formalism or Markov Chain models).

We showed numerical results obtained applying the proposed scheme to some models available in the literature, namely the Niederer-Hunter-Smith model [30], the model of Land and coworkers [27] and the mean-field model that we proposed in [43]. We considered a test case of an isotonic twitch and a test case where a contracted fiber shortens as a consequence of a gradual decrease of the external load, similarly to what happens in the ejection phase of the heart cycle. In both the test cases and for each force generation model here considered, the proposed scheme successfully removed the nonphysical numerical oscillations present in the solutions obtained with the standard segregated scheme, giving results in full agreement with the results of the monolithic scheme.

In order to quantify the numerical error introduced by the proposed stabilization term, we then performed a convergence analysis with respect to $\Delta t$. The results showed that, coherently with our analysis, the proposed scheme is convergent with order one. Moreover, the numerical error is only slightly larger than the error obtained with the monolithic scheme.

Finally, we applied the proposed scheme to a three-dimensional electromechanical simulation of the left ventricle. The results showed that the non-stabilized segregated scheme yields numerical oscillations affecting several variables of interest. Conversely, by adding the stabilization term proposed in this paper, these oscillations are successfully removed.

## Acknowledgements

This project has received funding from the European Research Council (ERC) under the European Union's Horizon 2020 research and innovation programme (grant agreement No 740132, iHEART - An Integrated Heart Model for the simulation of the cardiac function, P.I. Prof. A. Quarteroni).

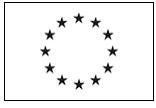
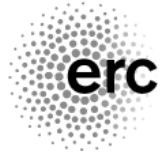

## A  Proofs of the stability results

We recall the following result (see e.g. [36, Property 1.13]):

**Proposition 5.** *Let $M \in \mathbb{C}^{n \times n}$ and $\varepsilon > 0$. There exists a norm $\|\cdot\|_{M,\varepsilon}$ that is compatible with the vector norm $|\cdot|$, i.e.*

$$|A\mathbf{x}| \leq \|A\|_{M,\varepsilon} |\mathbf{x}| \qquad \forall A \in \mathbb{C}^{n \times n}, \mathbf{x} \in \mathbb{C}^n,$$

*and satisfying*

$$\|M\|_{M,\varepsilon} \leq \rho(M) + \varepsilon.$$

We report below the proof of Props. 1 and 2.

*Proof of Prop. 1.* We denote $\eta := \sup_{k=0,\dots} |\boldsymbol{\eta}^{(k)}|$. Let us define the error $\mathbf{e}^{(k)} = \widetilde{\boldsymbol{\psi}}^{(k)} - \boldsymbol{\psi}^{(k)}$. Subtracting (27) from (28), we get $\mathbf{e}^{(0)} = \boldsymbol{\eta}^{(0)}$ and, for $k = 1, 2, \dots$:

$$\begin{aligned}
\mathbf{e}^{(k)} &= \Phi(\widetilde{\boldsymbol{\psi}}^{(k-1)}, t_k, \Delta t) - \Phi(\boldsymbol{\psi}^{(k-1)}, t_k, \Delta t) + \Delta t\, \boldsymbol{\eta}^{(k)} \\
&= \nabla_{\boldsymbol{\psi}} \Phi(\boldsymbol{\xi}^{(k)}, t_k, \Delta t)\, \mathbf{e}^{(k-1)} + \Delta t\, \boldsymbol{\eta}^{(k)},
\end{aligned} \tag{59}$$

for some $\boldsymbol{\xi}^{(k)}$ belonging to the line connecting $\boldsymbol{\psi}^{(k-1)}$ and $\widetilde{\boldsymbol{\psi}}^{(k-1)}$. Let us consider a finite time $T$ and let us set $\Delta t = T/N$. By Prop. 5, for any $N > 0$ and $k = 1, \dots, N$ there exists a compatible norm $\|\cdot\|_{N,k}$ such that $\|\nabla_{\boldsymbol{\psi}} \Phi(\boldsymbol{\xi}^{(k)}, t_k, \Delta t)\|_{N,k} \leq 1 + 2\alpha \Delta t$. It follows, for $k = 1, \dots, N$:

$$|\mathbf{e}^{(k)}| \leq \|\nabla_{\boldsymbol{\psi}} \Phi(\boldsymbol{\xi}^{(k)}, t_k, \Delta t)\|_{N,k} |\mathbf{e}^{(k-1)}| + \Delta t\, |\boldsymbol{\eta}^{(k)}| \leq (1 + 2\alpha \Delta t)|\mathbf{e}^{(k-1)}| + \Delta t\, |\boldsymbol{\eta}^{(k)}|. \tag{60}$$

Summing over $k$, we get, for $k = 1, \dots, N$:

$$|\mathbf{e}^{(k)}| \leq |\boldsymbol{\eta}^{(0)}| + \sum_{s=0}^{k-1} 2\alpha \Delta t |\mathbf{e}^{(k)}| + \sum_{s=0}^{k-1} \Delta t |\boldsymbol{\eta}^{(k)}| \tag{61}$$

From the Gronwall Lemma (see e.g. [36, Lemma 11.2]), it follows, for $k = 1, \dots, N$:

$$\begin{aligned}
|\mathbf{e}^{(k)}| &\leq \left(|\boldsymbol{\eta}^{(0)}| + \sum_{s=0}^{k-1} \Delta t |\boldsymbol{\eta}^{(k)}|\right) e^{2\alpha k \Delta t} \leq (1 + k\Delta t) e^{2\alpha k \Delta t} \eta \\
&\leq (1 + N\Delta t) e^{2\alpha N \Delta t} \eta = (1 + T) e^{2\alpha T} \eta.
\end{aligned} \tag{62}$$

This proves the zero-stability of the method, with $C = (1+T)e^{2\alpha T}$ (actually, the proof can be improved up to $C = (1 + T)e^{\alpha T}$). In order to prove the absolute stability, we notice that (59) entails

$$\mathbf{e}^{(k)} = \left[\prod_{s=1}^{k} \nabla_{\boldsymbol{\psi}}\Phi(\boldsymbol{\xi}^{(s)}, t_s, \Delta t)\right] \boldsymbol{\eta}^{(0)} + \Delta t \sum_{s=1}^{k} \left[\prod_{q=s+1}^{k} \nabla_{\boldsymbol{\psi}}\Phi(\boldsymbol{\xi}^{(q)}, t_q, \Delta t)\right] \boldsymbol{\eta}^{(k)}. \quad (63)$$

Let us consider a fixed $\Delta t$. By Prop. 5, for any $k$ there exists a compatible norm $\|\cdot\|_k$ such that $\|\nabla_{\boldsymbol{\psi}}\Phi(\boldsymbol{\xi}^{(k)}, t_k, \Delta t)\|_k \leq \rho_1$, for some $\rho_1$ such that $\rho_0 < \rho_1 < 1$. It follows

$$|\mathbf{e}^{(k)}| \leq \rho_1^k \eta + \Delta t(1 + \rho_1 + \rho_1^2 + \cdots + \rho_1^{k-1})\eta$$
$$= \rho_1^k \eta + \Delta t \frac{1 - \rho_1^k}{1 - \rho_1} \eta \xrightarrow[k \to +\infty]{} \Delta t \frac{\eta}{1 - \rho_1}. \quad (64)$$

Hence, the method is absolutely stable with constant $C = \Delta t/(1 - \rho_1)$. $\square$

*Proof of Prop. 2.* It is easily verified that, if $\mu_{0,0} \in [0, \mu_0^f/r]$, then the solution of (31) satisfies $\mu_0^{(k)} \in [0, \mu_0^f/r]$ for any $k$. Then, we can restrict our analysis to the values of $\boldsymbol{\psi} = (\mu_0, \mathbf{y}^T)^T$ such that $\mu_0 \in [0, \mu_0^f/r]$.

The transition map associated with the scheme (31) reads

$$\Phi((\mu_0, \mathbf{y}^T)^T, t_k, \Delta t) = \begin{pmatrix} \frac{\mu_0 + \mu_0^f \Delta t}{1 + r \Delta t} \\ C\left(\frac{\mu_0 + \mu_0^f \Delta t}{1 + r \Delta t}, \Delta t\right) \mathbf{y} + A\left(\frac{\mu_0 + \mu_0^f \Delta t}{1 + r \Delta t}, \Delta t\right)^{-1} \mathbf{h}^{(k)} \end{pmatrix} \quad (65)$$

By (H1), the determinant of $A(\mu_0, \Delta t)$ never vanishes. Hence, by the Cramer formula and by (H2), the map $\Phi$ is differentiable in its first argument. Its Jacobian matrix reads

$$\nabla_{\boldsymbol{\psi}}\Phi((\mu_0, \mathbf{y}^T)^T, t_k, \Delta t) = \begin{pmatrix} (1 + r \Delta t)^{-1} & \mathbf{0}^T \\ \boldsymbol{\zeta}(\mu_0, \mathbf{y}, \Delta t) & C\left(\frac{\mu_0 + \mu_0^f \Delta t}{1 + r \Delta t}, \Delta t\right) \end{pmatrix} \quad (66)$$

for a suitable function $\boldsymbol{\zeta}(\mu_0, \mathbf{y}, \Delta t)$. Hence, the spectrum of $\nabla_{\boldsymbol{\psi}}\Phi$ coincides with that of the matrix $C$, with the addition of the eigenvalue $(1 + r \Delta t)^{-1} \in (0, 1)$. Therefore, the spectral conditions (29) and (30) are equivalent to (35) and (36). Then the thesis follows from Prop. 1. $\square$

## B  Parameters of the distribution-moments model

In this appendix we list the values of the parameters of the distribution-moments model of Eq. (11) used in this paper. In order to derive realistic values for these parameters, we follow the calibration pipeline that we proposed in [42] for a generalized version of model (11).

Indeed, we have shown that the five parameters of the model considered in [42] can be chosen so to fit five experimentally measurable quantities, namely the isometric force, the fraction of attached binding sites in isometric conditions, the maximum shortening velocity, the curvature of the force-velocity relationship and the slope of the force-length response to fast transients.

We remark that the model considered in [42] has an additional parameter with respect to Eq. (11), that rules the curvature of the force-velocity relationship. With the model of Eq. (11), instead, the force-velocity relationship is always a straight line. Hence, we calibrate the four parameters of the model of Eq. (11) by fitting the experimentally measurable quantities listed above, except for the curvature of the force-velocity relationship. In this manner, by employing the same data reported in [42], we obtain the following values: $\mu_0^f = 114.4\,\text{s}^{-1}$, $\mu_1^f = 1.76\,\text{s}^{-1}$, $r = 520\,\text{s}^{-1}$ and $a_{\text{XB}} = 17.727\,\text{MPa}$.

## C  The NHS06 model

The NHS06 [30] model envisages five state variables:

$$\mathbf{r}(t) = \left([\text{Ca}^{2+}]_{\text{TRPN}}(t), z(t), Q_1(t), Q_2(t), Q_3(t)\right)^T,$$

where $[\text{Ca}^{2+}]_{\text{TRPN}}$ denotes the concentration of $[\text{Ca}^{2+}]$ bound to troponin, $z$ is fraction of available actin sites and $Q_1$, $Q_2$, $Q_3$ are the variables of the so-called *fading memory model*. The calcium-binding kinetics is ruled by the following ODE:

$$\frac{d[\text{Ca}^{2+}]_{\text{TRPN}}}{dt} = k_{\text{on}} \left([\text{Ca}^{2+}]_{\text{TRPN}}^{\max} - [\text{Ca}^{2+}]_{\text{TRPN}}\right) - k_{\text{off}} \left(1 - \frac{T_a}{\gamma_{\text{trpn}} T_{\text{ref}}}\right) [\text{Ca}^{2+}]_{\text{TRPN}},$$

where $k_{\text{on}}$ and $k_{\text{off}}$ are the association and dissociation rates of calcium, $T_{\text{ref}}$ is the reference tension and $\gamma_{\text{trpn}}$ rules the tension dependence. The fraction of available binding sites evolves according to

$$\dot{z} = \alpha_0 \left(\frac{[\text{Ca}^{2+}]_{\text{TRPN}}}{[\text{Ca}^{2+}]_{\text{TRPN}}^{50}(\lambda)}\right)^{n_H} (1 - z) - \alpha_{r1} z - \alpha_{r2} \frac{z^{n_{\text{rel}}}}{z^{n_{\text{rel}}} + K_z^{n_{\text{rel}}}},$$

where $[\text{Ca}^{2+}]_{\text{TRPN}}^{50}$, namely the calcium concentration corresponding to half-activation, is defined through

$$[\text{Ca}^{2+}]_{\text{TRPN}}^{50}(\lambda) = [\text{Ca}^{2+}]_{\text{TRPN}}^{\max} \frac{\text{Ca}_{50}(\lambda)}{\text{Ca}_{50}(\lambda) - \frac{k_{\text{off}}}{k_{\text{on}}}\left(1 - \frac{1+\beta_0\lambda}{2\gamma_{\text{trpn}}}\right)}, \quad (67)$$

$$\text{Ca}_{50}(\lambda) = \text{Ca}_{50}^{\text{ref}}(1 + \beta_1 \lambda),$$

and where $\alpha_0$, $\alpha_{r1}$ and $\alpha_{r2}$ are the rates of activation and relaxation (divided into slow and fast), $n_H$ and $n_{\text{rel}}$ are the nonlinear activation and relaxation coefficients and $\beta_1$ is a coefficient ruling the length-dependence.

The strain-dependence of the force generation mechanism is described through the following fading-memory:

$$\dot{Q}_i = A_i \dot{\lambda} - \alpha_i Q_i, \qquad i = 1, 2, 3, \qquad (68)$$

where $\alpha_i$ are the exponential rate constants and $A_i$ the associated weighting coefficients (for $i = 1, 2, 3$). Finally, the active tension predicted by this model is given by

$$T_a = T_{\text{ref}} (1 + \beta_0 \lambda) \frac{z}{z_{\max}(\lambda)} K\left(\sum_{i=1}^{3} Q_i\right),$$

where $\beta_0$ is a parameter ruling the length-dependence and where we have defined the maximum activation level associated with the strain $\lambda$ as

$$z_{\max}(\lambda) = \frac{\alpha_0 \left(\dfrac{[\text{Ca}^{2+}]_{\text{TRPN}}^{\max}}{[\text{Ca}^{2+}]_{\text{TRPN}}^{50}(\lambda)}\right)^{n_H} - K_2}{\alpha_0 \left(\dfrac{[\text{Ca}^{2+}]_{\text{TRPN}}^{\max}}{[\text{Ca}^{2+}]_{\text{TRPN}}^{50}(\lambda)}\right)^{n_H} + \alpha_{r1} + K_1}, \qquad (69)$$

$K_1$ and $K_2$ being suitable constants. The dependence of the generated force on the phenomenological variables of the fading-memory model is defined through the following function:

$$K(Q) = \begin{cases} \dfrac{aQ + 1}{1 - Q}, & \text{if } Q \leq 0, \\ \dfrac{(2 + a)Q + 1}{1 + Q}, & \text{if } Q > 0, \end{cases} \qquad (70)$$

where the constant $a$ controls the curvature of the force-velocity relationship and corresponds to the homonym parameter of the Hill model (see [18]).

The derivative of the function $K$, needed to define the active stiffness according to Eq. (47), is given by

$$K'(Q) = \frac{1 + a}{(1 + |Q|)^2}.$$

For the numerical discretization of the NHS06 model we consider a fully implicit scheme.

## D   The L17 model

In the L17 model [27], $\text{Ca}_{\text{TRPN}}(t)$ denotes the fraction of troponin units with calcium bound. Its evolution is described by

$$\frac{d\text{Ca}_{\text{TRPN}}}{dt} = k_{\text{TRPN}} \left( \left(\frac{[\text{Ca}^{2+}]_i}{[\text{Ca}^{2+}]_i^{50}(\lambda)}\right)^{n_{\text{TRPN}}} (1 - \text{Ca}_{\text{TRPN}}) - \text{Ca}_{\text{TRPN}} \right), \qquad (71)$$

where $k_{\text{TRPN}}$ is the kinetic constant and $n_{\text{TRPN}}$ is the cooperativity coefficient. The half-activating calcium concentration depends on the strain $\lambda$ through the following relationship:

$$[\text{Ca}^{2+}]_i^{50}(\lambda) = [\text{Ca}^{2+}]_i^{50,\text{ref}} + \beta_1 \min(\lambda, 0.2), \qquad (72)$$

where $\beta_1$ is a coefficient ruling the length-dependent calcium sensitivity.

The binding sites are split into four groups: blocked, unblocked (but without any XB attached), with weakly-bound XB (i.e. in pre-powerstroke configuration) and with strongly-bound XB (i.e. in post-powerstroke configuration). The fractions of sites in each of these groups are respectively denoted by the variables $B(t)$, $U(t)$, $W(t)$ and $S(t)$. Their evolution is described by the following system of ODEs:

$$\begin{cases} \dot{B} = k_b \, (\text{Ca}_{\text{TRPN}})^{-n_{\text{Tm}}/2} \, U - k_u \, (\text{Ca}_{\text{TRPN}})^{n_{\text{Tm}}/2} \, B \\ \dot{W} = k_{uw} U - (k_{wu} + \gamma_{wu}(\zeta_w) + k_{ws})W \\ \dot{S} = k_{ws} W - (k_{su} + \gamma_{su}(\zeta_s))S, \end{cases} \qquad (73)$$

while $U(t)$ can be obtained as $U(t) = 1 - B(t) - S(t) - W(t)$. In the above equations, $k_b$ and $k_u$ are the binding/unbinding rates, $n_{\text{Tm}}$ represents a cooperativity coefficient and $k_{\alpha\beta}$ (with $\alpha, \beta \in \{u, w, s\}$) are the transition rates among the XB configurations. Moreover, $\gamma_{su}$ and $\gamma_{wu}$ represent distortion-dependent unbinding rates, that depend on the variables $\zeta_w$ and $\zeta_s$, denoting the average distortion of XBs in the $W$ and $S$ groups, respectively. The time evolution of these variables is described by the following distortion-decay model:

$$\begin{cases} \dot{\zeta}_w = A_w \dot{\lambda} - c_w \zeta_w, \\ \dot{\zeta}_s = A_s \dot{\lambda} - c_s \zeta_s, \end{cases} \quad (74)$$

where $A_w$ and $A_s$ are the magnitude of the instantaneous response to distortion, while $c_w$ and $c_s$ are the decay rates of distortion. The distortion-dependent unbinding rate are defined as

$$\gamma_{wu}(\zeta_w) = \gamma_w |\zeta_w|, \qquad \gamma_{su}(\zeta_s) = \begin{cases} -\gamma_s(1 + \zeta_s), & \text{if } 1 + \zeta_s < 0 \\ \gamma_s \zeta_s, & \text{if } 1 + \zeta_s > 1 \\ 0, & \text{otherwise}. \end{cases}$$

Hence, the state vector contains the following six variables:

$$\mathbf{r}(t) = (\mathrm{Ca}_{\mathrm{TRPN}}(t), B(t), W(t), S(t), \zeta_w(t), \zeta_s(t))^T.$$

Finally, the generated active force is given by

$$T_{\mathrm{a}} = h(\lambda) \frac{T_{\mathrm{ref}}}{r_s} \left[ (1 + \zeta_s) S + \zeta_w W \right], \quad (75)$$

having defined the following functions:

$$h(\lambda) = \max(0, h_1(\min(1 + \lambda, 1.2))),$$
$$h_1(\psi) = 1 + \beta_0(\psi + \min(\psi, 0.87) - 1.87),$$

and where $\beta_0$ represents a length-dependence coefficient.

For the numerical discretization of the L17 model we consider two schemes. The first one is a fully implicit scheme. The second one is a semimplicit scheme, where we implicitly treat all the dependences, except for the variables $\mathrm{Ca}_{\mathrm{TRPN}}$, $\zeta_w$ and $\zeta_s$ in Eq. (73), that are explicitly treated.

## E  The RDQ20-MF model

In the RDQ20-MF model (corresponding to the model referred to as *MF-ODE model* in [43]), the thin filament activation is modeled by considering a triplet of consecutive regulatory units. Each regulatory unit is composed by troponin, that can be either unbound ($\mathcal{U}$) or bound ($\mathcal{B}$) to calcium, and by tropomyosin, that can be either in non-permissive ($\mathcal{N}$) or in permissive ($\mathcal{P}$) state. We denote by $\pi^{\alpha\beta\delta,\eta}$, where $\alpha, \beta, \delta \in \{\mathcal{N}, \mathcal{P}\}$ and $\eta \in \{\mathcal{U}, \mathcal{B}\}$ the probability that the tropomyosin units of the triplet are in $\alpha$, $\beta$ and $\delta$ state, respectively, and that the central troponin unit is in state $\eta$. The time evolution of these variables is described by the following ODE:

$$\begin{aligned}
\frac{d}{dt} \pi^{\alpha\beta\delta,\eta} &= \widetilde{k}_T^{\overline{\alpha}\alpha|\circ\cdot\beta,\circ} \pi^{\overline{\alpha}\beta\delta,\eta} - \widetilde{k}_T^{\alpha\overline{\alpha}|\circ\cdot\beta,\circ} \pi^{\alpha\beta\delta,\eta} \\
&+ k_T^{\overline{\beta}\beta|\alpha\cdot\delta,\eta} \pi^{\alpha\overline{\beta}\delta,\eta} - k_T^{\beta\overline{\beta}|\alpha\cdot\delta,\eta} \pi^{\alpha\beta\delta,\eta} \\
&+ \widetilde{k}_T^{\overline{\delta}\delta|\beta\cdot\circ,\circ} \pi^{\alpha\beta\overline{\delta},\eta} - \widetilde{k}_T^{\delta\overline{\delta}|\beta\cdot\circ,\circ} \pi^{\alpha\beta\delta,\eta} \\
&+ k_C^{\overline{\eta}\eta|\beta} \pi^{\alpha\beta\delta,\overline{\eta}} - k_C^{\eta\overline{\eta}|\beta} \pi^{\alpha\beta\delta,\eta},
\end{aligned} \quad (76)$$

where we use the notation $\overline{\mathcal{N}} = \mathcal{P}$, $\overline{\mathcal{P}} = \mathcal{N}$, $\overline{\mathcal{U}} = \mathcal{B}$ and $\overline{\mathcal{B}} = \mathcal{U}$ to denote opposite states, and where we set

$$\widetilde{k}_T^{\overline{\alpha}\alpha|\circ\,\cdot\,\beta,\circ} = \frac{\sum_{\xi,\zeta} k_T^{\overline{\alpha}\alpha|\xi\,\cdot\,\beta,\zeta} \pi^{\xi\overline{\alpha}\beta,\zeta}}{\sum_{\xi,\zeta} \pi^{\xi\overline{\alpha}\beta,\zeta}}, \qquad \widetilde{k}_T^{\overline{\delta}\delta|\beta\,\cdot\,\circ,\circ} = \frac{\sum_{\xi,\zeta} k_T^{\overline{\delta}\delta|\beta\,\cdot\,\xi,\zeta} \pi^{\beta\overline{\delta}\xi,\zeta}}{\sum_{\xi,\zeta} \pi^{\beta\overline{\delta}\xi,\zeta}}.$$

In the above equations, we denote by $k_C^{\eta\overline{\eta}|\beta}$ the troponin transition rate from state $\eta \in \{\mathcal{U}, \mathcal{B}\}$ to $\overline{\eta}$, when the corresponding tropomyosin unit is in state $\beta \in \{\mathcal{N}, \mathcal{P}\}$. Similarly, we denote by $k_T^{\overline{\beta}\beta|\alpha\,\cdot\,\delta,\eta}$ the tropomyosin transition rate from the state $\beta \in \{\mathcal{N}, \mathcal{P}\}$ to $\overline{\beta}$, when the adjacent tropomyosin units are in state $\alpha$ and $\delta$ and the corresponding troponin unit is in state $\eta \in \{\mathcal{U}, \mathcal{B}\}$.

The XB cycling dynamics is described within the formalism of the H57 model, under the assumption that the total attachment-detachment rate is independent of the myosin distortion, so that the distribution-moments reduction can be carried out. The population of binding sites is split into two groups, according to the state of the corresponding tropomyosin units. Hence, we denote by $\mu_{\mathcal{N}}^p$ and by $\mu_{\mathcal{P}}^p$ the $p$-th order moments of the population associated, respectively, with non-permissive and permissive tropomyosin units. The time evolution of these moments is ruled by the following system of ODEs:

$$\begin{cases} \dot{\mu}_{\mathcal{P}}^0 = -\left(r_0 + \alpha|\dot{\lambda}| + \widetilde{k}_T^{\mathcal{P}\mathcal{N}}\right)\mu_{\mathcal{P}}^0 + \widetilde{k}_T^{\mathcal{N}\mathcal{P}}\mu_{\mathcal{N}}^0 + P\,\mu_{f_{\mathcal{P}}}^0 \\ \dot{\mu}_{\mathcal{N}}^0 = -\left(r_0 + \alpha|\dot{\lambda}| + \widetilde{k}_T^{\mathcal{N}\mathcal{P}}\right)\mu_{\mathcal{N}}^0 + \widetilde{k}_T^{\mathcal{P}\mathcal{N}}\mu_{\mathcal{P}}^0 \\ \dot{\mu}_{\mathcal{P}}^1 = -\left(r_0 + \alpha|\dot{\lambda}| + \widetilde{k}_T^{\mathcal{P}\mathcal{N}}\right)\mu_{\mathcal{P}}^1 + \widetilde{k}_T^{\mathcal{N}\mathcal{P}}\mu_{\mathcal{N}}^1 + P\,\mu_{f_{\mathcal{P}}}^1 + \dfrac{d\lambda}{dt}\mu_{\mathcal{P}}^0 \\ \dot{\mu}_{\mathcal{N}}^1 = -\left(r_0 + \alpha|\dot{\lambda}| + \widetilde{k}_T^{\mathcal{N}\mathcal{P}}\right)\mu_{\mathcal{N}}^1 + \widetilde{k}_T^{\mathcal{P}\mathcal{N}}\mu_{\mathcal{P}}^1 + \dfrac{d\lambda}{dt}\mu_{\mathcal{N}}^0, \end{cases} \qquad (77)$$

where we have defined the permissivity $P$ (i.e. the fraction of units in state $\mathcal{P}$) as

$$P = \sum_{\alpha,\delta \in \{\mathcal{N},\mathcal{P}\}} \sum_{\eta \in \{\mathcal{U},\mathcal{B}\}} \pi^{\alpha\mathcal{P}\delta,\eta},$$

and where the following terms model the fluxes among the two populations:

$$\widetilde{k}_T^{\mathcal{N}\mathcal{P}} = \frac{\sum_{\alpha,\delta,\eta} k_T^{\mathcal{N}\mathcal{P}|\alpha\,\cdot\,\delta,\eta}\pi^{\alpha\mathcal{N}\delta,\eta}}{1-P}, \quad \widetilde{k}_T^{\mathcal{P}\mathcal{N}} = \frac{\sum_{\alpha,\delta,\eta} k_T^{\mathcal{P}\mathcal{N}|\alpha\,\cdot\,\delta,\eta}\pi^{\alpha\mathcal{P}\delta,\eta}}{P}.$$

In the above equations, $r_0$ denotes the total attachment/detachment rate in isometric conditions and $\alpha$ denotes the velocity-dependent detachment rate coefficient. Hence, the state of this model contains 20 variables (16 variables of the type of $\pi^{\alpha\beta\delta,\eta}$ and the four distribution-moments):

$$\mathbf{r}(t) = \left(\pi^{\mathcal{N}\mathcal{N}\mathcal{N},\mathcal{U}}, \pi^{\mathcal{P}\mathcal{N}\mathcal{N},\mathcal{U}}, \ldots, \pi^{\mathcal{P}\mathcal{P}\mathcal{P},\mathcal{B}}, \mu_{\mathcal{N}}^0, \mu_{\mathcal{P}}^0, \mu_{\mathcal{N}}^1, \mu_{\mathcal{P}}^1\right)^T.$$

Finally, the active tension predicted by the model is given by

$$T_{\mathrm{a}} = a_{\mathrm{XB}}\,\chi_{\mathrm{so}}(\lambda)\left[\mu_{\mathcal{P}}^1 + \mu_{\mathcal{N}}^1\right], \qquad (78)$$

where the single-overlap ratio $\chi_{\mathrm{so}}$ denotes the fraction of the AF filament in the single-overlap zone (the definition of the function $\chi_{\mathrm{so}}$ can be found in [43]). We remark that an implementation of the RDQ20-MF model is freely available in the following online repository:

<https://github.com/FrancescoRegazzoni/cardiac-activation>

For the numerical discretization of the RDQ20-MF model we consider a semimplicit scheme, as described in [43]. In this scheme, we employ an explicit treatment for the variables describing the thin-filament activation (i.e. the variables in the form of $\pi^{\alpha\beta\delta,\eta}$), while we implicitly treat the four moments $\mu_{\mathcal{N}}^0$, $\mu_{\mathcal{P}}^0$, $\mu_{\mathcal{N}}^1$ and $\mu_{\mathcal{P}}^1$.